\documentclass{amsart}
\usepackage{ifthen}
\usepackage{subfigure}
\usepackage{amsmath} 
\usepackage{amssymb} 
\usepackage{euscript}
\usepackage[all]{xy} 
\usepackage{varioref}
\usepackage{graphicx}
\usepackage{verbatim}
\usepackage{lscape}
\usepackage{xspace}
\usepackage{boxedminipage}
\usepackage{ifthen}
\newboolean{pdfOutput}
\setboolean{pdfOutput}{false}
\newboolean{abridged}
\setboolean{abridged}{false}
\ifthenelse{\boolean{pdfOutput}}
{
}
{
\usepackage{epsfig}
}
\newcommand{\drawfile}[1]
{
  \ifthenelse{\boolean{pdfOutput}}
  {
  \includegraphics{{#1}.pdf}%
  }
  {
  \epsfig{file={#1}.eps}%
  }
}
\newtheorem{theorem}{Theorem} 
\newtheorem{corollary}{Corollary}
\newtheorem{lemma}{Lemma} 
\newtheorem{proposition}{Proposition}

\theoremstyle{remark}
\newtheorem{definition}{Definition} 
\newtheorem{example}{Example}
\newtheorem{remark}{Remark}

\newcounter{remarkCounter}
\newlength{\setBracketHeight}

\newcommand{\LieDer}{\ensuremath{\EuScript L}}
\newcommand{\hook}{\ensuremath{\mathbin{ \hbox{\vrule height1.4pt
        width4pt depth-1pt \vrule height4pt width0.4pt depth-1pt}}}}

\newcommand{\pd}[2]{\ensuremath{\frac{\partial{#1}}{\partial{#2}}}}
\newcommand{\R}[1]{\ensuremath{\mathbb{R}^{#1}}}

\newcommand{\Gr}[2]{\ensuremath{\operatorname{Gr}\left({#1},{#2}\right)}}

\newcommand{\SO}[1]{\operatorname{SO}\left({#1}\right)}

\newcommand{\GL}[1]{\operatorname{GL}\left({#1}\right)}

\newcommand{\gl}[1]{\mathfrak{gl}\left({#1}\right)}
\newcommand{\SL}[1]{\operatorname{SL}\left({#1}\right)}

\newcommand{\PGL}[1]{\mathbb{P}\GL{#1}}
\newcommand{\slLie}[1]{\mathfrak{sl}\left({#1}\right)}

\newcommand{\Sym}[2]{\ensuremath{\operatorname{Sym}^{#1}\left({#2}\right)}}
\newcommand{\Lm}[2]{\ensuremath{\Lambda^{#1} \left ( {#2} \right )}}
\newcommand{\nForms}[2]{\ensuremath{\Omega^{#1} \left ( {#2} \right
    )}}

\DeclareMathOperator{\Ad}{Ad}

 \DeclareMathOperator{\tr}{tr}

 \DeclareMathOperator{\ad}{ad}
\DeclareMathOperator{\Aut}{Aut}

\newcommand{\Proj}[1]{\mathbb{P}^{#1}}
\newcommand{\ucProj}[1]{\widetilde{\mathbb{P}}^{#1}}
\newcommand{\Aff}[1]{\mathbb{A}^{#1}}

\newcommand{\OO}[1]{
  \ensuremath{
    \mathcal{O}
    \ifthenelse{\equal{#1}{0}}
      {}
      {\left({#1}\right)}
  }
}
\newcommand{\OOp}[2]{
  \ensuremath{
    \mathcal{O}
    \ifthenelse{\equal{#1}{0}}
      {}
      {\left({#1}\right)}
    \ifthenelse{\equal{#2}{1}}
      {}
      {^{\oplus{#2}}}
  }
}

\newcommand{\extraSection}[2]
{
\ifthenelse{\boolean{abridged}}
  {
  }
  {
    \section{#1} 
    \begin{center}
    \emph{This section will not be referred to
    subsequently, and may be skipped.}
    \end{center}
    \par{#2}
  }
}
\newcommand{\extraSubsection}[2]
{
\ifthenelse{\boolean{abridged}}
  {
  }
  {
    \subsection{#1}
    \begin{center}
    \emph{This subsection will not be referred to
    subsequently, and may be skipped.}
    \end{center}
    \par{#2}
  }
}
\newcommand{\extraStuff}[1]
{
\ifthenelse{\boolean{abridged}}
  {
  }
  {
    {#1}
  }
}
\newcommand{\F}[1]{\ensuremath{F{#1}}}
\newcommand{\Killing}{\ensuremath{\mathbb{B}}}

\begin{document}
\title{Complete projective connections}
\author{Benjamin McKay}
\address{University College Cork \\ Cork, Ireland} 
\email{b.mckay@ucc.ie}
\date{\today} 
\thanks{MSC {53B10}}
\thanks{Thanks (1) to two reviewers, for pointing out Kuiper's work,
and for correcting my misapprehensions of elementary facts
about Ricci curvature, (2) to David Wraith for explaining patiently 
the folklore of Ricci curvature, (3) to Andreas \v{C}ap for discussions on
left invariant Cartan geometries, and (4) to the Erwin Schr\"odinger
Institute for hospitality.}
\begin{abstract}
  The first examples of complete projective
  connections are uncovered: on surfaces, normal projective connections
  whose geodesics are all closed and embedded are complete.
  On manifolds of any dimension, normal projective connections induced from complete
  affine connections with slowly decaying positive Ricci curvature are complete.
\end{abstract}
\maketitle
\tableofcontents
\section{Introduction}
This article is a step toward global
analysis of Cartan geometries; a new avenue
of research, where almost nothing is known.
Completeness of projective connections is subtle, even on
compact manifolds, and there seems to be
no easy way to decide whether a projective
connection is complete. In my recent 
work \cite{McKay:2004}, I discovered
that complete complex projective connections
are flat, and I decided to look for complete
real projective connections which are not flat.
I was surprised to find that none were known;
in this article you will find the first
examples.\footnote{More recently I have discovered
that Tanaka \cite{Tanaka:1957} p. 21 announced in a remark
that he could prove that Einstein
metrics of positive Ricci curvature on compact manifolds are projectively
complete, although the proof did not appear. He actually says
negative Ricci curvature, but must clearly mean positive
Ricci curvature. It seems
likely that the intended proof is the same as mine.
I have also discovered that 
Blumenthal \cite{Blumenthal:1987} proved
that completeness of various
Cartan geometries is preserved under submersions, but 
this did not generate examples other than 
the standard Hopf fibrations.}

Definitions are presented later, but for the
moment recall that every Riemannian
manifold has a distinguished projective connection,
and that this imposes on every geodesic a natural
choice of parameterization, well defined modulo
projective transformations. 
This parameterization
is \emph{not} the arc length parameterization,
in many examples, but it is unchanged
if we change metric, as long as we keep the
same geodesics. 
For example, in real projective
spaces $\Proj{n}$, the projective parameterization
is the obvious parameterization: the geodesics are 
projective lines. But therefore in
affine space, thought of as an affine chart
of $\Proj{n}$ (which gives it the usual straight
lines as geodesics), following the projective 
parameterization can run 
us off to infinity in finite time. This is the
bizarre incompleteness of affine space. Worse:
since a flat torus is a quotient of affine space, it
is also incomplete! We wrap around
a geodesic infinitely often in finite time. Indeed, the only complete
examples known (before the results below) were
the sphere and projective space (with standard
metrics). 

The concept of completeness of Cartan geometries
is tricky to define, raised explicitly for the first time by
Ehresmann \cite{Ehresmann:1936} (also see 
Ehresmann \cite{Ehresmann:1938,Ehresmann:1951},
Kobayashi \cite{Kobayashi:1954}, 
Kobayashi \& Nagano \cite{KobayashiNagano:1964},
Clifton \cite{Clifton:1966}, Bates \cite{Bates:2004}),
and plays a central role in Sharpe's book \cite{Sharpe:1997},
but is also clearly visible beneath the surface in numerous works of
Cartan. Roughly speaking, completeness concerns
the ability to compare a geometry to some notion of
flat geometry, by rolling along curves. There were no
examples of complete projective connections
except for the sphere and projective space
(which are both flat) until now:
\begin{theorem}
Every normal projective connection on a surface, all of whose
geodesics are closed embedded curves, is complete.
\end{theorem}
\begin{theorem}
Every complete torsion-free affine connection with positive
Ricci curvature decaying slower than quadratically induces
a complete normal projective connection.
\end{theorem}
The first theorem is more exciting, since it is purely global and
depends directly on the projective connection.
\begin{example}
The projective connection on the product $S^n \times S^n$
of round spheres is projectively complete for $n > 1$. 
However, it is projectively incomplete for $n=1$, since
$S^1 \times S^1$ is the torus. But $S^1 \times S^1 \subset S^n \times S^n$
is a totally geodesic submanifold: a totally geodesic submanifold can have a different
projective parameterization from the projective parameterization
associated to the ambient manifold.
\end{example}
Analysis of projective connections is 
difficult because already on the simplest example,
projective space, the automorphism group
is not compact, a kind of inherent slipperiness.
A sort of antithesis of completeness is known
as \emph{projective hyperbolicity}; see
Kobayashi \cite{Kobayashi:1979} and Wu \cite{Wu:1981}
for examples of projective hyperbolicity.

{\'E}lie Cartan \cite{Cartan:70} introduced the notion of projective
connection; Kobayashi \& Nagano \cite{KobayashiNagano:1964},
Gunning \cite{Gunning:1978} and Borel \cite{Borel:2001} 
provide a contemporary review; we will use the definitions
of Kobayashi \& Nagano. This article may be difficult to follow
without the article of Kobayashi \& Nagano in hand. 
\section{The flat example: projective space}
First, let us consider projective space $\Proj{n}=\left(\R{n+1}\backslash 0\right)/\R{\times}$.  
Projective
space is glued together out of affine charts, and the
transition functions are affine transformations, so preserve straight
lines, i.e. geodesics.  The geodesic-preserving transformations of
projective space are precisely the projective linear transformations,
forming the group $\PGL{n+1,\R{}}$ (a well-known result in geometry due to
David Hilbert \cite{Hilbert:1999}).

We will think of $\Proj{n}$ as the space of tuples
\[
\begin{pmatrix}
  x^0 \\
  \vdots \\
  x^n \\
\end{pmatrix}
\]
of numbers, not all zero, modulo rescaling.  Write the corresponding
point of $\Proj{n}$ as
\[
\begin{bmatrix}
  x^0 \\
  \vdots \\
  x^n \\
\end{bmatrix}.
\] 
$\Proj{n}$ is acted on transitively by the group $G=\PGL{n+1,\R{}}$ of
projective linear transformations, i.e.  linear transformations of the
$x$ variables modulo rescaling. We will write $[g]$ for the element of
$\PGL{n+1,\R{}}$ determined by an element $g \in \GL{n+1,\R{}}.$ The stabilizer
of the point
\[
\begin{bmatrix}
  1 \\
  0 \\
  \vdots \\
  0 \\
\end{bmatrix}
\]
is the group $G_0$ consisting of $[g]$ where $g$ is a matrix
of the form
\[
[g] =
\begin{bmatrix}
  g^0_0 & g^0_j \\
  0 & g^i_j
\end{bmatrix}
\]
where $i,j=1,\dots,n.$ The Lie algebra of $\PGL{n+1, \R{}}$ is just
$\slLie{n+1, \R{}},$ so consists of the matrices of the form
\[
\begin{pmatrix}
  A^0_0 & A^0_j \\
  A^i_0 & A^i_j
\end{pmatrix}
\]
with $A^0_0 + A^i_i = 0.$ We define the Maurer--Cartan 1-form $\Omega
\in \nForms{1}{\PGL{n+1, \R{}}} \otimes \slLie{n+1, \R{}}$ by \( \Omega = g^{-1}
\, dg.  \) This form satisfies $d \Omega = - \Omega \wedge \Omega.$
Splitting into components, we calculate
\begin{align*}
  d \Omega^i_0 &= - \left ( \Omega^i_j + \delta^i_j \Omega^k_k \right
  )
  \wedge \Omega^j_0 \\
  d \Omega^i_j &= - \Omega^i_k \wedge \Omega^k_j
  + \Omega^0_j \wedge \Omega^i_0 \\
  d \Omega^0_i &= \left( \Omega^j_i + \delta^j_i \Omega^k_k \right )
  \wedge \Omega^0_j
\end{align*}
Following Cartan \cite{Cartan:1992} we let $\omega^i = \Omega^i_0, \gamma^i_j = \Omega^i_j + \delta^i_j
\Omega^k_k,$ and $\Omega_i = \Omega^0_i$ then we find
\begin{align*}
  d \omega^i &= - \gamma^i_j \wedge \omega^j \\
  d \gamma^i_j &= - \gamma^i_k \wedge \gamma^k_j
  + \left ( \omega_j \delta^i_k + \omega_k \delta^i_j \right ) \wedge \omega^k \\
  d \omega_i &= \gamma^j_i \wedge \omega_j.
\end{align*}

The group $G_0$ is a semidirect product: each element
factors into two elements of the form
\[
\begin{bmatrix}
  1 & 0 \\
  0 & g
\end{bmatrix}
\begin{bmatrix}
  1 & \lambda \\
  0 & 1
\end{bmatrix}. 
\]
It will be helpful later to see how each of these factors acts on our
differential forms. This is not difficult, since the form
$\Omega=g^{-1} \, dg$ satisfies
\[
r_{g_0}^* \Omega = \Ad_{g_0}^{-1} \Omega,
\]
for $g_0 \in G_0$
We leave to the reader to calculate that if we write $g$ for the
matrix
\[
\begin{bmatrix}
  1 & 0 \\
  0 & g
\end{bmatrix}
\]
and $\lambda$ for the matrix
\[
\begin{bmatrix}
  1 & \lambda \\
  0 & 1
\end{bmatrix}
\]
then
\begin{align*}
  r_g^* \omega^i &= \left(g^{-1}\right)^i_j \omega^j \\
  r_g^* \gamma^i_j &= \left(g^{-1}\right)^i_k \gamma^k_l g^l_j \\
  r_g^* \omega_i &= \omega_j g^j_i \\
  r_{\lambda}^* \omega^i &= \omega^i \\
  r_{\lambda}^* \gamma^i_j &= \gamma^i_j + \left ( \lambda_j
    \delta^i_k + \lambda_k \delta^i_j \right )
  \omega^k \\
  r_{\lambda}^* \omega_i &= \omega_i - \lambda_j \gamma^j_i -
  \lambda_i \lambda_j \omega^j.
\end{align*}

We can reconsider projective geometry in terms of bundles. For any
manifold $M$ of dimension $n$, let $\F{M}$ (called the \emph{frame bundle}
of $M$) be the set of all isomorphisms of tangent spaces
of $M$ with $\R{n}$.
The group $G=\PGL{n+1, \R{}}$ acts transitively on $\Proj{n},$ and also on the
frame bundle $\F{\Proj{n}}.$ 
The stabilizer of a point of $\Proj{n}$ is
$G_0 \subset G$; the stabilizer of a frame at a point is
the subgroup $G_1 \subset G_0$ consisting of matrices of the form
\[
\begin{bmatrix}
  1 & \lambda \\
  0 & 1
\end{bmatrix}.
\]
It is easy to show that $\Proj{n}$ has
tangent spaces $T_P \Proj{n} = P^* \otimes \left(\R{n+1}/P\right)$.

We can identify
\[
\xymatrix{
  \F{\Proj{n}} \ar@{<->}[r] \ar[d] & \PGL{n+1,\R{}}/G_1 \ar[d] \\
  \Proj{n} \ar@{<->}[r] & \PGL{n+1}/G_0.  }
\]
We have another bundle over $\Proj{n}$, $\PGL{n+1, \R{}}$ itself, which we
can put on the top at the right side.  We will build a corresponding
bundle on the left side.

Consider the geodesics of projective space.  These are the projective
lines. If we think of projective space as the space of lines through 0
in a vector space, its geodesics correspond to 2-planes in that vector
space.  Thus the space of geodesics is
$\Gr{2}{n+1}=\PGL{n+1, \R{}}/G_2$ where $G_2$
consists of the matrices of the form
\[
[g] =
\begin{bmatrix}
  g^0_0 & g^0_1 & g^0_J \\
  g^1_0 & g^1_1 & g^1_J \\
  0 & 0 & g^I_J
\end{bmatrix}
\]
where $I,J=2,\dots,n$.  Above the space of geodesics is the space of
pointed geodesics, which is the space of choices of a 2-plane in our
vector space with a line in that 2-plane, so it is
$\PGL{n+1}/G_+$ where $G_+ \subset
G_0$ consists of matrices of the form
\[
[g] =
\begin{bmatrix}
  g^0_0 & g^0_1 & g^0_J \\
  0 & g^1_1 & g^1_J \\
  0 & 0 & g^I_J
\end{bmatrix}.
\]
Write $\mathfrak{g}$ and $\mathfrak{g}_0$ for the Lie algebras
of $G$ and $G_0$, etc.
\section{Structure equations of projective connections}
Given a right principal $G_0$-bundle $E \to M$,
write $r_g : E \to E$ for the right $G_0$-action
of $g \in G_0$. We will refer to $G_0$ as the
\emph{structure group} of the principal bundle.
\begin{definition}
A \emph{projective connection} on an $n$-manifold $M$ is a 
choice of principal right $G_0$-bundle
$E \to M$ together with 
a 1-form $\Omega \in \nForms{1}{E} \otimes \mathfrak{g}$
so that
\begin{enumerate}
\item at each point $e \in E$, $\Omega_e : T_e E \to \mathfrak{g}$
is a linear isomorphism
\item $r_g^* \Omega=\Ad^{-1}_g \Omega$, and
\item for any $A \in \mathfrak{g}$, writing $\vec{A}$ 
for the unique vector field satisfying $\vec{A} \hook \Omega = A$,
we require further that
\[
e^{\vec{A}}=r_{e^{A}}
\]
whenever $A \in \mathfrak{g}_0$ (the left
hand side is the flow of a vector field).
\end{enumerate}
\end{definition}
Write
\[
\Omega =
\begin{pmatrix}
\Omega^0_0 & \Omega^0_j \\
\Omega^i_0 & \Omega^i_j
\end{pmatrix}
\]
with $\Omega^0_0 + \Omega^i_i = 0$.
We follow Cartan and define 1-forms $\omega^i,\gamma^i_j,\omega_i$
($i,j,k,l=1,\dots,n$), linearly independent,
by the equations
\begin{align*}
\omega^i &= \Omega^i_0 \\
\gamma^i_j &= \Omega^i_j - \delta^i_j \Omega^0_0 \\
\omega_j &= \Omega^0_j.
\end{align*}
(This is just a change of basis from the $\Omega^{\bullet}_{\bullet}$
1-forms.) Since $\Omega$ is a 1-form valued in 
$\mathfrak{g}$, we can think of $\omega^{\bullet}$ as
$\Omega \mod \mathfrak{g}_0$.
\begin{lemma}
There are uniquely determined functions $K^i_{jk},K_{ijk},K^i_{jkl}$
(called the \emph{curvature functions})
so that the structure equations of Cartan 
in table~\vref{tbl:Struc} are satisfied.
\end{lemma}
\begin{proof}
This requires elementary applications
of Cartan's lemma; see Kobayashi \& Nagano~\cite{KobayashiNagano:1964} for proof.
\end{proof}
\begin{table}
\begin{align*}
  \nabla \omega^i &= d \omega^i + \gamma^i_j \wedge \omega^j \\
  &= \frac{1}{2} K^i_{kl} \omega^k \wedge \omega^l \\
  \nabla \gamma^i_j &= d \gamma^i_j + \gamma^i_k \wedge \gamma^k_j -
  \left( \omega_j \delta^i_k + \omega_k \delta^i_j
  \right) \wedge \omega^k  \\
  &= \frac{1}{2} K^i_{jkl} \omega^k \wedge \omega^l \\
  \nabla \omega_i &= d \omega_i - \gamma^j_i \wedge \omega_j \\
  &= \frac{1}{2} K_{ikl} \omega^k \wedge \omega^l \\
  0 &= K^i_{jk} + K^i_{kj} \\
  0 &= K^i_{jkl} + K^i_{jlk} \\
  0 &= K_{ikl} + K_{ilk}. \\
\end{align*}
\caption{The structure equations of a projective connection}\label{tbl:Struc}
\end{table}
\section{Elementary global aspects of projective connections}
\begin{definition}
A projective connection is called \emph{flat} if the curvature functions vanish.
\end{definition}
\begin{example}
The model of a projective connection is the one on
$\Proj{n}$ given by taking $E=G$, $\Omega=g^{-1} \, dg$ 
the left invariant Maurer--Cartan 1-form, and the
map $g \in \PGL{n+1,\R{}} \to g \left[e_0\right] \in \Proj{n}$.
The model is flat. 
\end{example}
\begin{lemma}
A projective connection is flat just when it is locally
(i.e. on open subsets of $M$) isomorphic to the 
model.
\end{lemma}
\begin{proof}
Clearly local isomorphism implies flatness. Start with a flat projective
connection.Take the exterior differential system 
$\Omega-g^{-1} \, dg=0$ on $E \times \PGL{n+1,\R{}}$.
(See Bryant et al. \cite{BCGGG:1991} for more on
exterior differential systems, Cauchy characteristics,
and integral manifolds.) It satisfies the conditions of the Frobenius
theorem just when the curvature functions
vanish. The $\mathfrak{g}$ orbits
are Cauchy characteristics, so maximal
connected integal manifolds are unions
of these orbits. The
group $G_0$ has finitely many path components,
and the union of finitely many integral manifolds
is an integral manifold, so each connected
integral manifold is contained in a unique
$G_0$-invariant and $\mathfrak{g}$-invariant
integral manifold.
Integral manifolds are the graphs
of local isomorphisms, which descend to maps
$M \to \Proj{n}$ by $G_0$-equivariance.
\end{proof}
\begin{example}
The sphere $S^n$ has a 2-1 covering map $S^n \to \Proj{n}$.
Pulling back the bundle from the model, and the 1-form
$g^{-1} \, dg$ from the model, we find a flat projective
connection on $S^n$.
\end{example}
\begin{definition}[Ehresmann \cite{Ehresmann:1951}]
A projective connection is \emph{complete}
when the vector fields $\vec{A}$ are all
complete (i.e their flows are defined
for all time).
\end{definition}
\begin{example}
The vector fields $\vec{A}$ on the model
generate the right action of $\PGL{n+1,\R{}}$
on itself, and therefore the model is complete.
\end{example}
\begin{example}
A covering space of a complete projective
connection is complete, because the relevant
vector fields are pullbacks under covering maps.
Therefore $S^n$ is complete.
\end{example}
\begin{definition}
An \emph{isomorphism} of projective connections
$E_j \to M_j$ $(j=0,1)$ is a $G$-equivariant diffeomorphism
$\Phi : E_0 \to E_1$ preserving the projective
connection forms: $\Phi^* \Omega_1 = \Omega_0$.
An \emph{infinitesimal symmetry} of 
a projective connection $E \to M$ is a vector
field $X$ on $E$ commuting with the $G_0$-action,
and satisfying $\LieDer_X \Omega=0$.
\end{definition}
\begin{lemma}\label{lemma:CompleteSymmetries}
If a projective connection is complete then
every infinitesimal symmetry is a complete vector field. 
\end{lemma}
\begin{proof}
(Essentially the same as Bates \cite{Bates:2004}.)
The vector fields $\vec{A}$
commute with $X$, so they permute the flow lines
of $X$ around in all directions. Therefore the
time for which the flow of $X$ is defined
is locally constant. But then it cannot diminish
as we move along a flow line. Therefore the
flow of $X$ is defined for all time.
\end{proof}
\begin{lemma}
A projective connection $E \to M$ is flat
just when the infinitesimal symmetries act locally 
transitively on $E$, complete and flat just when
the automorphism group is transitive on $E$.
\end{lemma}
\begin{proof}
This forces invariance of the curvature,
which therefore must be constant. 
The curvature is equivariant under the
$G_0$ action, so lives in a $G_0$-representation.
One can easily see that there are no nonzero
$G_0$-invariant vectors in that representation.
\end{proof}
\begin{theorem}\label{thm:pushPull}
Every flat projective connection $E \to M$
is obtained by taking the universal covering space
$\tilde{M} \to M$, mapping $\tilde{M} \to \Proj{n}$
by a local diffeomorphism,
pulling back the model projective connection,
and taking the quotient projective connection
via some morphism $\pi_1(M) \to \PGL{n+1,\R{}}$.
\end{theorem}
\begin{remark}
The map $\tilde{M} \to \Proj{n}$ is called
the \emph{developing map}.
\end{remark}
\begin{proof}
Without loss of generality, assume that 
$M$ is connected.
Put the exterior differential system $\Omega=g^{-1} \, dg$
on the manifold $E \times \PGL{n+1,\R{}}$.
By the Frobenius theorem, the manifold
is foliated by leaves (maximal connected
integral manifolds). Because the system
is invariant under left action of
$\PGL{n+1,\R{}}$ on itself, this
action permutes leaves.
Define vector fields $\vec{A}$
on $E \times \PGL{n+1,\R{}}$
by adding the one from 
$E$ with the one (by the same name)
from $\PGL{n+1,\R{}}$.
The flow of $\vec{A}$ on $\PGL{n+1,\R{}}$
is defined for all time, so
the vector field $\vec{A}$
on $E \times \PGL{n+1,\R{}}$
has flow through a point 
$(e,g)$ defined for as long 
as the flow is defined down on $E$.
These vector fields $\vec{A}$ are Cauchy
characteristics, so the leaves are
invariant under their flows.

The group $G_0$ has
finitely many components, so
the $G_0$ orbit of a leaf is a finite
union of leaves. Let $\Lambda_0$
and $\Lambda_1$ be $G_0$-orbits of leaves, 
containing points
$\left(e_i,g_i\right) \in \Lambda_i$.
After replacing these points by
other points obtained through $G_0$ action, 
we can draw a path
from $e_0$ to $e_1$ in $E$, consisting
of finitely many flows of $\vec{A}$ vector 
fields, so such a path lifts to
our leaf. Therefore 
$\Lambda_1$ must contain a point
$\left(e_0,g_0'\right)$. Therefore there
is a $G_0$-orbit $\Lambda$ of a leaf,
unique up to $\PGL{n+1,\R{}}$ action.

The inclusion
$\Lambda \subset E \times \PGL{n+1,\R{}}$
defines two local diffeomorphisms
$\Lambda \to E$ and $\Lambda \to \PGL{n+1,\R{}}$,
both $\vec{A}$ and $G_0$ equivariant. Consider
the first of these. Let $F$ be a fiber
of $\Lambda \to E$ over some point $e \in E$.
Define local coordinates on $E$ by
inverting the map
$A \in \mathfrak{g} \mapsto e^{\vec{A}} e \in E$
near $A=0$. This map is only defined
near $A=0$, and is a diffeomorphism
in some neighborhood, say $U$, of $0$. Then map 
\[
U \times F \to \Lambda 
\] 
by $(A,f) \mapsto e^{\vec{A}}f$, clearly a local
diffeomorphism. Therefore $\Lambda \to E$
is a covering map, and $G_0$-equivariant,
so descends to a covering map $\tilde{M}=\Lambda/G_0 \to M=E/G_0$.
Thus $\Lambda \to \tilde{M}$ is the pullback
bundle of $E \to M$.

The map $\Lambda \to \PGL{n+1,\R{}}$ is $G_0$-equivariant, 
so descends to a map $\tilde{M} \to \Proj{n}$.
By definition, on $\Lambda$ we have $\Omega=g^{-1} \, dg$,
so this map is pullback of projective connections.
\end{proof}
\begin{corollary}
A complete flat projective connection $E \to M$ with $\dim M \ge 2$ 
is isomorphic to a quotient $S^n/\Gamma$ with $\Gamma \subset \SO{n+1}$
a finite group.
\end{corollary}
\begin{proof}
The universal covering space has to be the same as for $\Proj{n}$,
so $S^n$. But then the quotient has to be by a discrete 
group $\Gamma$ with morphism $\Gamma \to \SL{n+1,\R{}}$. Since $S^n$
is compact, it can only act as covering space of 
compact spaces, so with a finite group $\Gamma$
of deck transformations. Every finite subgroup
of $\SL{n+1,\R{}}$ preserves a positive definite
inner product on $\R{n+1}$, so $\Gamma$ sits
in a conjugate of $\SO{n+1}$. 
\end{proof}
\section{Classification of projective connections on curves}
\begin{lemma}
Every projective connection $E \to M$ on a curve $M$ is flat.
\end{lemma}
\begin{proof}
Curvature is a semibasic 2-form, but $M$ has
only one dimension.
\end{proof}
Consider $\Proj{1}=\Aff{} \cup \infty$, $\Aff{}=\R{}$ and let $\Aff{+}$ be the 
positive real numbers. We draw the universal cover $\ucProj{1} \to \Proj{1}$
as
\[
\xymatrix{
\dots \ar@{-}[rrrrr] & \ar@{}[]^{0} & \ar@{}[]^{\infty} & \ar@{}[]^{0} & \ar@{}[]^{\infty} & \dots \\
}
\]
Pull back the standard flat projective connection on $\Proj{1}$ to a projective
connection on $\ucProj{1}$. 
The automorphism group 
$\Aut \ucProj{1}$ of that projective 
connection on $\ucProj{1}$ is the obvious 
central extension of $\PGL{2,\R{}}$ by $\mathbb{Z}$,
where $\mathbb{Z}$ acts by translating zeros
to zeros in this picture; write this action
as $n \in \mathbb{Z} : x \mapsto x+\infty_n$. 
To be more concrete, we can split up each
$2 \times 2$ matrix $g$ into $g=qr$, the usual
$QR$-factorization from linear algebra.
For $g \in \SL{2,\R{}}$, $g=qr$ with
\[
q =
\begin{pmatrix}
\cos \phi & - \sin \phi \\
\sin \phi & \cos \phi
\end{pmatrix},
\
r
=
\begin{pmatrix}
a & b \\
0 & 1/a
\end{pmatrix}
\]
with $a > 0$. Think of $a,b,\phi$ as local
coordinates on $\SL{2,\R{}}$. (In terms of
these coordinates, group operations are unbearably
complicated.) The $a,b,\phi$ are clearly global coordinates
on the universal covering group of $\SL{2,\R{}}$,
which is $\Aut \ucProj{1}$. Moreover,
$r \in G_0$, so $\phi$ is a global
coordinate function on $\ucProj{1}=\Aut \ucProj{1}/G_0$,
quotienting out the right $G_0$ action.

In terms of the standard affine chart, identifying
\[
x \in \R{} \to 
\begin{bmatrix}
x \\
1
\end{bmatrix}
\in \Proj{1},
\]
$x=\cot(\phi)$ maps $\ucProj{1} \to \Proj{1}$.
Then $x=0$ lies at $\phi=\pi/2$, and $x=\infty$ lies at $\phi=0$.
The affine (left) group action
on $\Proj{1}$, $x \mapsto mx+b$, lifts to a unique action
on $\ucProj{1}$ fixing all $\infty$'s. This is
just the use of elements of $G_0$ on the left
instead of the right. We see this by direct calculation:
$\phi=0$ or $\phi=\pi$ just when $q$ commutes with $r$.
Thus $qr$ represents the same point of $\ucProj{1}$
as $q$. 

\begin{theorem}[Kuiper \cite{Kuiper:1954}, Gorinov \cite{Gorinov:2006}]
The projective connections on a closed connected curve
(modulo isomorphism) are:
\par{}\noindent{}%
\begin{tabular}{ll} 
elliptic & 
$\ucProj{1}/\left(\phi \mapsto \phi + \theta \right)$
\\
parabolic & (1) $ \ \Aff{}/(x \mapsto x+1)$ or \\
          & (2) $ \ \ucProj{1}/(x \mapsto x+\infty_n+1)$ \\
hyperbolic & (1) $ \ \Aff{+}/(x \mapsto rx)$ (some $r>1$) or \\
	& (2) $ \ \ucProj{1}/(x \mapsto rx+\infty_n)$
\end{tabular}
\par{}\noindent{}%
where we can assume that $n$ is an arbitrary positive integer,
the angle $\theta$ can be any nonzero real number,
and $n \ne 0$.
In particular, the numbers $r, \theta, n$ are
invariants of the projective connection. 
The projective connections on an open connected curve 
(modulo isomorphism) are the pullbacks to the 
following open subsets of $\ucProj{1}$:
\par{}\noindent{}%
\begin{tabular}{ll}
elliptic &
	\xymatrix{& \ucProj{1}}
\\
parabolic 
	&
	\xymatrix{
		\text{\emph{(1)}} \ar@{}[r] & 
		\ar@{}[]^{-\infty} \ar@{-}[r] & 
		\dots \ar@{-}[r] &  
		\ar@{}[]^{\infty}
	} \\
	& 
	\xymatrix{
		\text{\emph{(2)}} \ar@{}[r] & 
		\ar@{}[]^{-\infty} \ar@{-}[r] & 
		\dots \ar@{-}[r] & 
		\ar@{}[]^{\infty} & 
		\ar@{-}[ll] \dots
	}
\\
hyperbolic &
\xymatrix{\ar@{}[r] & \ar@{}[]^{0} \ar@{-}[r] & \dots \ar@{-}[r] &  \ar@{}[]^{\infty}
}
\end{tabular}
\par{}\noindent{}%
Any number of copies of the affine line may
be contained in these open parabolic and hyperbolic curves
in the $\dots$ in the middle, and infinitely many must
appear in any $\dots$ at each end.
The hyperbolic curves are those with an
invariant metric $\left(dz/z\right)^2$ (and
consequently an invariant affine connection: the Levi-Civita
connection), defined except at $z=\infty$ and $z=0$. The parabolic are those with no invariant metric
but have an invariant affine connection, defined except at $z=\infty$, and the
elliptic are those which have neither.
\end{theorem}
\begin{proof} We will only outline the proof.
The technique is to use the developing map, i.e. identify the universal
cover locally with $\Proj{1}$, and thereby globally with
an open interval of $\ucProj{1}$, following theorem~\vref{thm:pushPull}.
If this open interval has an endpoint, we can slide it along
by automorphisms of $\ucProj{1}$, and put it where we like.
If it has two endpoints, we have to be more careful: we can
put one of them where we like, say at some $\infty$,
but then if the other one winds up landing at another
$\infty$ in the process, it is impossible to move it without
moving the first one. On the other hand, if the second
end point does not land on an $\infty$, we can slide
it along by affine transformations to land on a $0$.
For an open curve, this finishes the story.
Consider a closed curve.
With this normalization completed, a certain
subgroup of automorphisms is still available fixing
the (0, 1 or 2) endpoints. With this, we can normalize
the monodromy around the closed curve, which is
an element of $\PGL{2,\R{}}$. But the monodromy
must act without fixed points in the interior of 
$\tilde{C}$, while fixing all of the endpoints.
This allows us to classify the possible monodromy
elements up to conjugation by automorphisms.
\end{proof}
\begin{theorem}
A projective connection on a curve is complete just when
its universal cover is identified with $\ucProj{1}$ by
the developing map, i.e. either elliptic or
closed parabolic of type (2) or closed hyperbolic of type (2).
\end{theorem}
\begin{proof}
Let $E \to C$ be a projective connection on a curve.
Following theorem~\vref{thm:pushPull}, the universal cover $\tilde{C}$ is
mapped locally diffeomorphically to $\ucProj{1}$, and therefore
is a connected open subset. The bundle $E \to C$ lifts
to a bundle $\tilde{E} \to \tilde{C}$.
Completeness is invariant under covering maps, so
$E \to C$ is complete just when $\tilde{E} \to \tilde{C}$ is.
Clearly $\tilde{E}$ is an open subset of the automorphism group 
of $\ucProj{1}$, under
\[
\xymatrix{
E \ar[d] & \tilde{E} \ar[l] \ar[d] \ar[r] & \Aut \ucProj{1} \ar[d] \ar[r]  & \PGL{2,\R{}} \ar[d] \\
C        & \tilde{C} \ar[l] \ar[r]        & \ucProj{1} \ar[r]              & \Proj{1}. 
}
\]
Completeness is just precisely invariance of that open subset
under the flow of all left invariant vector fields,
i.e. under left translation by the identity component,
i.e. the open subset being a union of components. 
But $\widetilde{\PGL{2,\R{}}}=\Aut \ucProj{1}$ has precisely two components, 
and they are interchanged by right $G_0$ action,
and $\tilde{E}$ is $G_0$ invariant.
Therefore completeness 
is just equality of $\tilde{E}$ and $\Aut \ucProj{1}$.
So completeness of $C$ is just 
completeness of $\tilde{C}$ which is just isomorphism
of $\tilde{C}$ with $\ucProj{1}$.
\end{proof}
\section{Geodesics}
\begin{definition}
Given an immersed curve $\iota : C \to M$
on a manifold $M$ with projective connection
$E \to M$, define the pullback bundle $\iota^* E$
\[
\xymatrix{
\iota^* E \ar[d] \ar[r] & E \ar[d] \\
C \ar[r] & M.
}
\]
Inside $\iota^* E$, $\omega$ has rank
1 and transforms as $r_g^* \omega = g^{-1} \omega$
under the structure group. Therefore there
is a subbundle $E_C \subset \iota^* E$ on
which $\omega^I=0$, a principal $G_+$-subbundle,
where $G_+ \subset G_0$ is the subgroup
of projective transformations
preserving the projective line through
$\left[e_0\right]$ and $\left[e_1\right]$
as well as fixing the point $\left[e_0\right]$.
Taking exterior derivative of the equations
$\omega^I=0$, we find $\gamma^I_1 = \kappa^I \omega^1$,
for some functions $\kappa^I : E_C \to \mathbb{R}$,
which descend to a section $\kappa$ of 
$\left(\iota^* TM/TC\right) \otimes \left(T^*C\right)^2$,
called the \emph{geodesic curvature} of $C$.
(See Cartan \cite{Cartan:1992}, {L}e\c cons sur la th\'eorie des espaces \`a connexion
              projective, pp. 91-111.)
\end{definition}
\begin{remark}
To prove that $\kappa^I$ descend to a section
of this bundle, or prove other similar statements,
the procedure is always the same as our proof
that $TM = E \times_{G_0} \R{n}$ in lemma~\vref{prop:TM}.
\end{remark}
\begin{definition}
We define a \emph{geodesic} to be a curve of vanishing
geodesic curvature. 
\end{definition}
Equivalently, geodesics are the curves on $M$
which are the projections to $M$ of the integral
manifolds $E_C \subset E$ of the differential
system $\omega^I=\gamma^I_1=0$. 
\begin{definition}
The \emph{geodesic flow} is the flow of
the vector field dual to $\omega^1$.
\end{definition}
The flow lines of geodesic flow are contained in these
integral manifolds; since these
flow lines are permuted by the action of 
$G_+$, the manifolds $E_C$ for $C$ a geodesic
are precisely the $G_+$-orbits of flow
lines of geodesic flow.
\begin{definition}
Given a connected immersed 
curve $\iota : C \subset M$, and a chosen
point $c_0 \in C$, we will \emph{roll}
an immersed curve $C$ onto $\Proj{n}$ to produce an
immersion $C \to \Proj{n}$ (called
its \emph{development}), as follows. Take the 
differential system $\Omega=g^{-1} \, dg$
on $E \times \PGL{n+1,\R{}}$, restrict
it to $E_C \times \PGL{n+1,\R{}}$.
The Frobenius theorem once again tells
us that $E_C \times \PGL{n+1,\R{}}$
is foliated by leaves (i.e. maximal
connected integral manifolds 
of that differential system).

Take the $G_+$ orbit of any leaf,
say $\Lambda \subset E_C \times \PGL{n+1,\R{}}$,
containing a point $(e,1)$ with $\pi(e)=c$,
and map $\Lambda \to \PGL{n+1,\R{}} \to \Proj{n}$
by the obvious maps. By $G_+$-equivariance,
this determines a map $C \to \Proj{n}$.
By $\PGL{n+1,\R{}}$ invariance, changing the choice of $\Lambda$
changes the map $C \to \Proj{n}$ by a projective automorphism.
\end{definition}

\begin{definition}
Take an immersed curve $\iota : C \to M$
and a projective connection $E \to M$.
Let $N_+ \subset G_+$ be the subgroup
of $G_+$ acting trivially on $\Proj{1} \subset \Proj{n}$.
Then $\bar{E}_C = E_C/N_+ \to C$
is a principal right $\bar{G}_0$-bundle,
where $\bar{G}_0=G_+/N_+$. 
Let $\bar{G}=G/N_+$. The 1-form
$\bar{\Omega}=\Omega \mod \mathfrak{n}_+ \in \nForms{1}{\bar{E}_C} \otimes \bar{\mathfrak{g}}$
is a projective connection on $\bar{E}_C \to C$.
Call this the \emph{induced} projective connection
on $C$.
\end{definition}

\begin{lemma}
$C$ is a geodesic just when
\begin{enumerate}
\item the development $C \to \Proj{n}$
maps $C \to \Proj{1} \subset \Proj{n}$, and 
\item
the induced projective connection
is the pullback via $C \to \Proj{1}$.
\end{enumerate}
\end{lemma}
\begin{proof}
This is immediate from the structure equations.
\end{proof}
It is elementary to prove:
\begin{theorem}[Kobayashi \cite{Kobayashi:1954}]
Let $M$ be a manifold with projective
connection.
Every immersed curve in projective space
is the development of a curve in $M$
just when $M$ is complete.
\end{theorem}

\begin{definition}
A curve $C$ in a manifold $M$ with projective connection 
is \emph{complete} when the induced projective connection $\bar{E}_C$
is complete.
\end{definition}
\begin{lemma}\label{lemma:geodCompleteness}
A projective connection is complete just when all
of its geodesics are complete.
\end{lemma}
\begin{proof}
The geodesic flow is tangent to all of the manifolds $E_C$ for
all geodesics, and projects under $E_C \to \bar{E}_C$
to the geodesic flow.
Therefore its completeness is identical to the
completeness of all geodesics. But 
$r_g^* \vec{A} = \overrightarrow\Ad_g A$ for
$g \in G_0$, permuting the vector fields dual
to all of the $\omega^i$. The vector fields
$\vec{A}$ for $A \in \mathfrak{g}_0$ are 
always complete, moving up the fibers. 
\end{proof}
\begin{remark}
It is unknown whether there are manifolds
of dimension greater than one 
with all geodesics parabolic and
closed, or with all geodesics elliptic and open.
\end{remark}
\section{Affine connections, projective connections, projective structures
and normal projective connections}
Just to clarify a few minor points in the literature,
we would like to explain the relations between
affine connections, projective connections, projective
structures and normal projective connections.
Kobayashi \& Nagano \cite{KobayashiNagano:1964}
explain how to relate projective structures,
torsion-free affine connections, and
normal projective connections, but they
don't explain the relation between 
arbitrary projective connections and
arbitrary affine connections.

Let $E \to M$ be a projective connection.
Clearly the bundle $E/\R{n*} \to M$ is a principal right $\GL{n,\R{}}$-bundle.

Let $\F{M} \to M$ be the bundle whose fiber over a 
point $m \in M$ consists in the linear isomorphisms
$u : T_m M \to \R{n}$. Make this into a principal
right $\GL{n,\R{}}$-bundle, by defining
$r_g u = g^{-1} u$ for $g \in \GL{n,\R{}}$.
For each $e \in E,$ $\Omega_e : T_e E \to \mathfrak{g}$
is onto, so if we write $\omega_e : T_e E \to \mathfrak{g}/\mathfrak{g}_0=\R{n}$
for the composition with the obvious projection, then
$\omega \in \nForms{1}{E} \otimes \R{n}$.
But $\omega=0$ on vertical vectors on $E$, so
for each point $e \in E$, $\omega_e$ determines
a linear isomorphism $\underline{\omega}_e : T_m M \to \R{n}$,
where $e \in E_m$.  
\begin{lemma}
The map $e \in E \to \underline{\omega}_e \in \F{M}$ descends to a 
$\GL{n,\R{}}$-equivariant
bundle isomorphism $E/\R{n*} \to \F{M}$. 
Let $\phi : \F{M} \to M$ be the bundle map.
Define 1-forms $\omega^i$ on $\F{M}$
by
\[
v \hook \omega = u(\phi'(u) v)
\]
for $v \in T_u \F{M}$. Then $E \to \F{M}$
pulls back $\omega^i$ to $\omega^i$
\end{lemma}
\begin{proof}
The $\GL{n,\R{}}$-equivariance is a calculation:
\[
\underline{\omega}\left(r_g e\right) = g^{-1} \underline{\omega}(e),
\]
from which the rest easily follows.
\end{proof}
\begin{lemma}
If $E \to M$ is a projective connection, then
$E \to E/\R{n*}=\F{M}$ is a trivial principal bundle right $\R{n*}$-bundle.
The $\GL{n,\R{}}$-equivariant sections of this bundle determine
affine connections on $M$ with the
given geodesics and given torsion $K^i_{jk}$. Changing the choice of
section alters the parameterization of the geodesics.
\end{lemma}
\begin{proof}
Existence of a global section $s$ is
elementary, for any principal $\R{n*}$-bundle, 
using local convex combinations.

Let $s : \F{M} \to E$ be a local
$\GL{n,\R{}}$-equivariant section.
Then $s^* \omega^i = \omega^i$,
clearly. Lets write
$\gamma^i_j$ still for the 1-forms
$s^* \gamma^i_j$, and $\omega_i$
for $s^* \omega_i$.

The vector fields
$\vec{A}$ on $E$, for $A \in \gl{n,\R{}}$,
project to the corresponding vector fields
$\vec{A}$ on $\F{M}$ given by the
right action of $\GL{n,\R{}}$. Therefore
$\vec{A} \hook \gamma=A$ on both
$E$ and $\F{M}$. Therefore $\vec{A} \hook \omega_i=0$
on both $E$ and $\F{M}$. 

The vectors dual to $\omega^i$ on $E$
project to nonzero vectors on $M$, and therefore
on $\F{M}$, because $\omega^i$ are semibasic.
Therefore $\omega^i,\gamma^i_j$ form a coframing
on $\F{M}$. So $\omega_i$ must be a combination
of them, and since $\vec{A} \hook \omega_i=0$
for $A \in \gl{n,\R{}}$, we must have $\omega_i=a_{ij} \omega^j$
for some functions $a_{ij} : \F{M} \to \R{}$.

Pick a geodesic $C \subset M$. Then
the tangent spaces of $E_C$ are cut 
out by the equations $\omega^I=\gamma^I_1=0$.
These equations are expressed in semibasic
1-forms, so the integral manifolds will project
to integral manifolds of the same system 
on $\F{M}$. The projective parameterizations of a 
geodesic are those given by the geodesic flow through
points of $E_C$, and therefore
unless the section $s$ stays inside a region
where $\underline{\omega}^1$ is constant,
the parameterization will
not match a projective parameterization. 

Because $\vec{A} \hook \gamma = A$
for $A \in \gl{n,\R{}}$, $\gamma$
determines a unique connection on $\F{M}$,
with horizontal space $\gamma=0$.
We leave the reader to show that the
geodesics of the connection are 
the integral manifolds of the exterior
differential system $\omega^I=\gamma^I_1=0$.

Recall that the torsion of an affine
connection with connection 1-forms $\gamma^i_j$
is given by equivariant functions $T^i_{jk}$, where
\[
d \omega^i = - \gamma^i_j \wedge \omega^j + \frac{1}{2} T^i_{jk} \omega^j \wedge \omega^k,
\]
so torsion is $T^i_{jk}=K^i_{jk}$, same as for
the projective connection.
\end{proof}
The structure equations given by
the choice of some section $s$ are
\begin{align*}
\nabla_s \omega^i &= d \omega^i + \gamma^i_j \wedge \omega^j \\
&= \frac{1}{2} K^i_{kl} \omega^k \wedge \omega^l \\
\nabla_s \gamma^i_j &= d \gamma^i_j + \gamma^i_k \wedge \gamma^k_j \\
&= \left( \frac{1}{2} K^i_{jkl} - a_{jl} \delta^i_k - a_{kl} \delta^i_j \right) 
\omega^k \wedge \omega^l,
\end{align*}
relating the curvature of the affine connection
to the curvature of the projective connection:
\[
R^i_{jkl} = K^i_{jkl} - 2 \left( a_{jl} \delta^i_k + a_{kl} \delta^i_j \right).
\]
\begin{corollary}
The unparameterized geodesics of a projective connection are
the unparameterized geodesics of some affine connection.
\end{corollary}
\begin{lemma}
Given a projective connection, there is a torsion-free
projective connection (i.e. $K^i_{jk}=0$) with the same 
geodesics, with the same projective parameterizations.
\end{lemma}
\begin{proof}
Set $\tilde{\gamma}^i_j = \gamma^i_j + \frac{1}{2} K^i_{jk}$.
\end{proof}
Kobayashi \cite{Kobayashi:1995} and Cartan \cite{Cartan:70}
show that 
given a projective connection, there is a
unique normal projective connection with the
same unparameterized geodesics. By our result
above, there are
torsion-free connections with the same
parameterized geodesics as this normal
projective connection.
Two affine connections are said to be 
\emph{projectively equivalent} if they
have the same unparameterized geodesics,
and an equivalence class is called 
a \emph{projective structure}.
Kobayashi \cite{Kobayashi:1995}
and Kobayashi \& Nagano \cite{KobayashiNagano:1964} also
show that given any torsion-free affine connection,
there is a unique bundle $E \to \F{M}$,
whose sections are precisely the torsion-free affine connections,
and that this bundle bears a unique normal
projective connection with the given geodesics.
Therefore a normal projective connection
is essentially the same object as a projective
structure. 
\begin{lemma}[Weyl]
Two connections $\gamma,\tilde{\gamma}$ on $\F{M}$
have the same geodesics up to parameterization
just when 
\[
\tilde{\gamma}^i_j = \gamma^i_j + \left(
\lambda_j \delta^i_k + \lambda_k \delta^i_j \right)
+ a^i_{jk} \omega^k
\]
where $\lambda_j \omega^j$ is the pullback
to $\F{M}$ of a 1-form $\lambda$ on $M$,
and $a^i_{jk} \omega^j \wedge \omega^k \otimes \pd{}{\omega^i}$
is the pullback to $\F{M}$ of a section
of $\Lm{2}{T^*M}\otimes TM$.
\end{lemma}
\begin{proof}
Any two connection 1-forms have to agree on
the vertical vectors, so can only differ
by semibasic 1-forms:
\[
\tilde{\gamma}^i_j = \gamma^i_j + p^i_{jk} \omega^k.
\]
The equations of geodesics of $\tilde{\gamma}^i_j$
are $\omega^I=\tilde{\gamma}^I_1=0$, giving
$\omega^I=\gamma^I_1 + p^I_{11} \omega^k$. 
For these to be the same Frobenius exterior
differential systems, their leaves must have 
the same tangent spaces. The connections will
share a geodesic $C$ just when the submanifolds
of $\F{M}$:
\[
\xymatrix{
\F{M}_C \ar[d] \ar[r] & \F{M} \ar[d] \\
C \ar[r] & M
}
\]
are the same.
These are the leaves of the exterior
differential system, which satisfies
the conditions of the Frobenius theorem, 
i.e. the conormal bundle is spanned
precisely by the 1-forms in the
exterior differential system,
so the systems must be
identical. Therefore $p^I_{11}=0$, for all
$I>1$. Since indices can be freely permuted
in this argument, we have $p^i_{jj}=0$
whenever $i \ne j$. By $\GL{n,\R{}}$-equivariance
of connection 1-forms, $p^i_{jj}=0$ whenever $i \ne j$.
Check that the expression
\[
P \hook \omega^i = p^i_{jk} \omega^j \wedge \omega^k
\]
defines a section $P$ of $T^*M \otimes T^*M \otimes TM$,
by $\GL{n,\R{}}$-equivariance. Split into 
$D=A+S$ with $A$ antisymmetric and $S$ symmetric.
Then $S(v,v)$ must be a multiple of $v$ for
all vectors $v \in TM$, because $p^i_{jj}=0$.
Therefore $S(v,v)=\lambda(v)v$ for a unique
1-form $\lambda$. Moreover, 
\[
\tilde{\gamma}^i_j=\gamma^i_j + \left(\lambda_j \delta^i_k + \lambda_k \delta^i_j \right) \omega^k
+ a^i_{jk} \omega^k.
\]
Equivariance under $\GL{n,\R{}}$ ensures that 
$a^i_{jk}$ determines a section of the appropriate
bundle on $M$.
\end{proof}
\begin{corollary}[Weyl]
Two connections $\gamma,\tilde{\gamma}$ on $\F{M}$
have the same geodesics up to parameterization
and the same torsion just when 
\[
\tilde{\gamma}^i_j = \gamma^i_j + \left(
\lambda_j \delta^i_k + \lambda_k \delta^i_j \right)
\]
where $\lambda_j \omega^j$ is the pullback
to $\F{M}$ of a 1-form $\lambda$ on $M$.
\end{corollary}
\begin{theorem}
Let $\nabla$ be an affine connection on a manifold $M$, and 
let $\underline{\gamma}$ be the corresponding connection
1-form on $\F{M}$. 
Define the 1-forms $\omega^i$ as above on $\F{M}$.
Define $E$ to be the set of triples
$\left(m,u,\tilde{\gamma}\right)$ so that $m \in M$, $u : T_m M \to \R{n}$
is a linear isomorphism,  and $\tilde{\gamma}$ is the value at
$u \in \F{M}$ of a smooth connection 1-form with the same
torsion and geodesics as $\underline{\gamma}$. 
Identify $\lambda \in \R{n*}$ with the element
of $\left(\lambda_j \delta^i_k + \lambda_k \delta^i_j\right) \in \R{n*} \otimes \gl{n,\R{}}$.
Make $\R{n*}$ act
on $E$ on the right by $r_{\lambda}(m,u,\tilde{\gamma})=\left(m,u,\tilde{\gamma}+\lambda \cdot \omega\right)$.
Then $E$ has the structure of a smooth manifold,
and $\phi : (m,u,\tilde{\gamma}) \in E \to u \in \F{M}$ is a principal
right $\R{n*}$-bundle, and $\Phi : (m,u,\tilde{\gamma}) \in E \to m \in M$
is a principal right $G_0$-bundle. 
Let the vector fields $\vec{\lambda}$ and $\vec{A}$
on $E$ be the  generators of the right action.
Define 1-forms on $E$
by pulling back $\omega^i$, defining
$\gamma^i_j$ by the equation
\[
v \hook \gamma^i_j = \left(\phi'(m,u,\tilde{\gamma}) v \right) \hook \tilde{\gamma}^i_j.
\]
Then there are 1-forms $\omega_i$ on $E$
so that 
\[
\vec{\lambda} \hook \omega_i = \lambda_i \text{ and } \vec{A} \hook \omega_i = 0
\]
for $A \in \gl{n,\R{}}$. We can pick these
1-forms uniquely if we require that
\[
d \gamma^i_j = -\gamma^i_k \wedge \gamma^k_j + 
d \gamma^i_j =- \gamma^i_k \wedge \gamma^k_j +
  \left( \omega_j \delta^i_k + \omega_k \delta^i_j
  \right) \wedge \omega^k  + \frac{1}{2} K^i_{jkl} \omega^k \wedge \omega^l
\]
for some functions $K^i_{jkl}$ and require further that $K^i_{jkl}+K^i_{jlk}=0$
and $K^i_{jil}=0$. Moreover, $\omega^i,\gamma^i_j,\omega_i$
constitute a projective connection $E \to M$ whose
geodesics are the geodesics of $\nabla$, and whose
torsion is the torsion of $\nabla$. This projective
connection is normal just when $\nabla$ is torsion-free.
\end{theorem}
\begin{proof}
It is obvious that $E \to \F{M}$ is a principal right 
$\R{n*}$-bundle. Define the
right action of $g \in \GL{n,\R{}}$ by
\[
r_g \left( m, u, \tilde{\Gamma} \right)
=
\left(
m,g^{-1}u, 
\Ad_g^{-1} \tilde{\Gamma}r_{g*}^{-1}
\right).
\]
Check that this fits together with the 
$\R{n*}$-action into a $G_0$-action.
For $g \in \GL{n,\R{}}$, check that
\begin{align*}
r_g^* \omega^{\bullet} &= g^{-1} \omega^{\bullet} \\
r_g^* \gamma &= \Ad_g^{-1} \gamma.
\end{align*}
Differentiate these equations to show that
for $A \in \gl{n,\R{}}$,
\begin{align*}
\LieDer_{\vec{A}} \omega^{\bullet} &= - A \omega^{\bullet} \\
\LieDer_{\vec{A}} \gamma &= - \Ad_A \gamma.
\end{align*}
Similarly, for $\lambda \in \R{n*}$,
\begin{align*}
r_{\lambda}^* \omega^{\bullet} &= \omega^{\bullet} \\
r_{\lambda}^* \gamma &= \gamma + \lambda \omega,
\end{align*}
so that
\begin{align*}
\LieDer_{\vec{\lambda}} \omega^{\bullet} &= 0 \\
\LieDer_{\vec{\lambda}} \gamma &= \lambda \cdot \omega^{\bullet}. 
\end{align*}
Since the $\vec{A}$ and $\vec{\lambda}$ exhaust
the vertical directions, Cartan's lemma
ensures that:
\begin{align*}
0 &= d \omega^i + \gamma^i_j \wedge \omega^j + \frac{1}{2} K^i_{kl} \omega^k \wedge \omega^l \\
0 &= d \gamma^i_j + \gamma^i_k \wedge \gamma^k_j -
\left( \omega_k \delta^i_j + \omega_j \delta^i_k \right) \wedge \omega^k
+
\frac{1}{2} K^i_{jkl} \omega^k \wedge \omega^l,
\end{align*}
for some functions $K^i_{kl}, K^i_{jkl}$ antisymmetric in $k,l$,
and some 1-forms $\omega_i$ with $\vec{\lambda} \hook \omega_i = \lambda_i$.

Taking the initial $\underline{\gamma}$ 1-form to determine
a section $\sigma(u)=\left(m,u,\underline{\gamma}\right)$ of $E \to \F{M}$,
calculate that
\begin{align*}
\sigma^* \omega^{\bullet} &= \omega^{\bullet} \\
\sigma^* \gamma &= \underline{\gamma},
\end{align*}
which implies that
\begin{align*}
0 &= \sigma^* \left( d \omega^i + \gamma^i_j \wedge \omega^j  + \frac{1}{2} K^i_{kl} \omega^k \wedge \omega^l \right) \\
 &= d \omega^i + \underline{\gamma}^i_j \wedge \omega^j + \frac{1}{2} K^i_{kl} \omega^k \wedge \omega^l, 
\end{align*}
so that along that section of $E$, we have
\[
K^i_{kl} = T^i_{kl}.
\]
Check that
\[
r_{\lambda}^* K^i_{kl} = K^i_{kl},
\]
so that $K^i_{kl}=T^i_{kl}$ everywhere on $E$.

There is still some freedom to pick these
$\omega_i$. We can change them to
\[
\tilde{\omega}_i = a_{ij} \omega^j
\]
without changing the equations we have
developed so far. This will alter the
expression $K^i_{jil}$, changing it to
\[
\tilde{K}^i_{jil} = K^i_{jil} + 2 a_{jl}.
\]
Therefore there is a unique choice of
$\omega_i$ 1-forms for which $K^i_{jil}=0$. 
Taking exterior derivative of the equations
so far, we find that the structure equations
of a projective connection are satisfied.
Moreover, the $G_0$-equivariance is assured
because the condition which determined
the $\omega_i$ was $G_0$-invariant.
\end{proof}
\begin{remark}
Notice the peculiar condition that $K^i_{jil}=0$.
This vanishing of what we might call \emph{``unsymmetrized 
Ricci curvature''} is required to specify
the choice of projective connection uniquely.
\end{remark}
\begin{corollary}
Any two projective connections on a manifold $M$, with the
same torsion and unparameterized geodesics 
as a given affine connection, differ
by a section of $T^*M \otimes T^*M$.
\end{corollary}
\begin{proof}
If $\omega_i$ is changed to $\tilde{\omega}_i$
in the above proof, then $\tilde{\omega}_i - \omega_i = a_{ij} \omega^j$,
and the structure equations tell us that
\begin{align*}
r_g^* a_{ij} &= a_{kl} g^i_k g^j_l \\
r_{\lambda}^* a_{ij} &= a_{ij},
\end{align*}
which ensures that $a_{ij}$ descends to a section
of that bundle.
\end{proof}
\section{Vector bundles and descent data}
\begin{definition}
Let $E \to M$ be a projective connection,
with connection 1-form $\Omega \in \nForms{E}\otimes \mathfrak{g}.$
We have a 1-form $\omega = \left(\omega^i\right) = \Omega \mod \mathfrak{g}_0
: T_e E \to \mathfrak{g}/\mathfrak{g}_0  = \R{n}.$ 
The linear map $\Omega : T_e E \to \mathfrak{g}$ 
is an isomorphism, identifying the vertical
directions with $\mathfrak{g}_0$. 
Therefore the $\omega^i$ are semibasic, i.e. vanish on the vertical
directions, being valued in $\R{n}=\mathfrak{g}/\mathfrak{g}_0$.
Moreover the 1-forms $\omega^i$ are a basis for the semibasic
1-forms for the map $E \to M$. At each point
$e \in E$, $\omega^i$ is therefore the pullback
via $\pi : E \to M$ of some 1-form from $\pi(e)$,
say $\underline{\omega}^i : T_{\pi(e)} M \to \R{}$,
so that 
\[
v \hook \omega = \pi(e)'v \hook \underline{\omega}^i
\] 
for all $v \in T_e E$. These $\underline{\omega}^i$
are \emph{not} sections of $T^*M$, but of $\pi^* T^*M$.
\end{definition}
\begin{lemma}
If $X$ is a vector field on $M$, define
functions $X^i : E \to \R{}$ by 
\[
X^i(e) = X \hook \underline{\omega}^i(e).
\]
Let
\[
X^{\bullet} = 
\begin{pmatrix}
X^1 \\
X^2 \\
\vdots \\
X^n
\end{pmatrix}
: E \to \R{n}.
\]
Then $r_g^* X^{\bullet} = g^{-1} X^{\bullet}$,
for $g \in G_0$.
\end{lemma}
\begin{proof}
Pick any vector field $Y$ on $E$ so that $\pi'(e)Y(e)=X(\pi(e))$.
We can do this by picking $Y$ locally, and making affine combinations
of local choices of $Y$. 
Therefore $r_{g*} Y - Y$ is vertical, for any $g \in G_0$.
The $G_0$-equivariance of $\Omega$
says that
\[
r_g^* \omega^{\bullet} = g^{-1} \omega^{\bullet},
\]
i.e.
\[
r_{g*} v \hook \omega^{\bullet}\left(r_g e\right) = g^{-1} \left( v \hook \omega^{\bullet}(e) \right).
\]
Therefore
\begin{align*}
r_g^* X^{\bullet}(e) 
&= X^{\bullet}\left(r_g e\right) \\
&= X \hook \underline{\omega}^{\bullet}\left(r_g e \right) \\
&= Y \hook \omega^{\bullet}\left(r_g e\right) \\
&= g^{-1} \left( \left( r_{g*}^{-1} Y \hook \omega^{\bullet}(e) \right) \right) \\
&= g^{-1} \left( Y(e) \hook \omega^{\bullet}(e) \right) \\
&= g^{-1} X^{\bullet}(e).
\end{align*}
\end{proof}
\begin{definition} Given a group $G_0$ and two spaces
on which $G_0$ acts on the right, say $X$ and $Y$, the
diagonal right $G_0$-action is the one given by $\left(x,y\right)g_0=\left(xg_0,yg_0\right)$.
Let $X \times_{G_0} Y$ be the quotient by the diagonal $G_0$-action.
\end{definition}
\begin{proposition}\label{prop:TM}
$TM = E \times_{G_0} \R{n}$
\end{proposition}
\begin{proof}
Given a vector field $X$ on $M$, 
the functions $X^{\bullet} : E \to \R{n}$
are $G_0$-equivariant, so form a section
of $E \times_{G_0} \R{n} \to M$.
Clearly $X=0$ just where that section vanishes,
so this is an injection of vector bundles of
equal rank, hence an isomorphism.
Alternately, given $G_0$-equivariant function $X^i$,
we define a section of $\pi^* TM \to E$ (where 
$\pi : E \to M$ is our projective connection
bundle), by $X \hook \underline{\omega}^{\bullet} = X^{\bullet}$.
But by $G_0$-equvariance of $X^{\bullet}$ and of $\omega^{\bullet}$,
this $X$ is $G_0$-invariant, so drops to a section
of $TM \to M$. 
\end{proof}
\begin{remark}
We will frequently state that various
equivariant expressions on various
principal bundles determine sections
of various vector bundles, and in each
case a proof along the lines of the above
applies, so we will omit those proofs. For example:
\end{remark}
\begin{corollary}
$T^*M = E \times_{G_0} \R{n*}$, etc.
\end{corollary}
\section{Positive Ricci curvature}
\begin{theorem}
Let $M$ have a complete affine connection
with Ricci curvature tensor $\operatorname{Ric}_{ij}=\frac{1}{2}\left(R^k_{ikj}+R^k_{jki}\right)$
positive definite, and bounded
from below along any geodesic by
\[
\operatorname{Ric}_{ij}(t) \ge \frac{c}{4t^{2-\epsilon}} \delta_{ij},
\]
for some constants $c,\epsilon > 0$
(these constants possibly dependent on the choice
of geodesic), for all sufficiently
large values of $t$, where $t$ is 
the natural affine parameter along
a geodesic determined by the affine
connection. Then the induced projective
connection on $M$ is complete.
\end{theorem}
\begin{proof}
Pick coordinates $x^i$ along the geodesic,
which we can assume (by parallel transport
of a frame) are constant on parallel transported
vectors. There will be infinitely many zeroes
to solutions of the equation
\[
y''+\operatorname{Ric}_{ij} \frac{dx^i}{dt} \frac{dx^j}{dt} y=0
\]
where $x^i$ are local coordinates along the 
geodesic, by comparison to the equation
\[
y''+\frac{c}{t^{2-\epsilon}}y=0,
\]
whose solutions are constant linear combinations of
\[
\sqrt{t} J\left(\frac{1}{\epsilon},\frac{\sqrt{c} t^{\epsilon/2}}{\epsilon}\right),
\sqrt{t} Y\left(\frac{1}{\epsilon},\frac{\sqrt{c}t^{\epsilon/2}}{\epsilon }\right),
\]
where $J$ and $Y$ are Bessel functions.
Kobayashi \& Sasaki \cite{Kobayashi/Sasaki:1979}
proved that the projective
parameterizations are precisely the ratios $u=(ay_1+b)/(cy_0+d)$,
for $y=y_0,y=y_1$ any two linearly independent solutions
of this equation, and $a,b,c,d$ any constants with
$ad-bc \ne 0$. Therefore the projective parameterization
wraps around infinitely often, ensuring completeness
of each geodesic by comparison to the classification
of projective connections on curves. Completeness
of every geodesic ensures completeness of the ambient
projective connection by lemma~\vref{lemma:geodCompleteness}.
\end{proof}
\begin{remark}
Every example known of such an affine
connection is found on a compact manifold. Kobayashi \& Nagano
\cite{KobayashiNagano:1964} wonder whether
projective completeness implies
compactness, which would imply a strengthened
Bonnet--Myers--Cheng theorem,
i.e. that slower than quadratic Ricci curvature
decay of a complete affine connection would
force compactness. For Riemannian manifolds,
there does not appear to be in the literature
any proof that slower than quadratic Ricci curvature
decay of the Levi-Civita connection would
force compactness,
but is not difficult to prove using results of
David Wraith \cite{Wraith:2005}.
\end{remark}
\begin{remark}
This theorem gives
rise to many examples, but depends on Ricci
curvature, which is not an invariant of a 
projective connection. We would like to 
find a criterion for completeness which 
can be checked in many examples, and
which is projectively invariant. Even for the projective
connections of Riemannian
manifolds, it is unclear to what extent Ricci curvature
bounds are really required. In fact they
are not: we will provide examples of surfaces
with projectively
complete Riemannian metrics, whose curvature takes
on both positive and negative values.
\end{remark}
\begin{example} The Killing form metric on a compact semisimple Lie group 
has positive Ricci curvature (Milnor \cite{Milnor:1976}), and therefore
is projectively complete.
\end{example}
\begin{remark} This theorem is similar to Tanaka \cite{Tanaka:1957} p. 21,
but he uses the opposite sign convention for Ricci curvature (so that 
the sphere has negative Ricci curvature for him), and his result requires
invariance of the Ricci curvature under parallel transport.
\end{remark}
\section{Left invariant projective connections on Lie groups}
\begin{theorem}
The isomorphism classes of left invariant projective connections on a Lie group $H$
of dimension $n$ with Lie algebra $\mathfrak{h}$
are invariantly identified with the linear maps $\mathfrak{h} \to \slLie{n+1,\R{}}$
which are transverse to $\mathfrak{g}_0$, modulo the Lie algebra automorphisms
of $\slLie{n+1,\R{}}$ fixing $\mathfrak{g}_0$, and the Lie algebra automorphisms
of $\mathfrak{h}$.
\end{theorem}
\begin{proof}
Suppose that $E \to H$ is a left
invariant projective connection, i.e. that the left action of $H$ on
itself lifts to a left action on $E$, preserving a projective connection.
The actions of $H$ and $G_0$ must commute.
Pick a point $e \in E$ and map $(h,g) \in H \times G_0 \to r_{g_0} he \in E$.
This map is clearly a diffeomorphism, so henceforth identify $E= H \times G_0$. 
Consider the projective connection $\Omega$. Lets write the
left invariant Maurer--Cartan 1-form on $G_0$ as $g_0^{-1} \, dg_0$,
and similarly write the
left invariant Maurer--Cartan 1-form on $H$ as $h^{-1} \, dh$.
Clearly $\Omega-g_0^{-1} \, dg_0$ vanishes on the fibers of 
$H \times G_0$, so $\Omega - g_0^{-1} \, dg_0$ is a multiple
of $h^{-1} \, dh$. Lets write it as 
\begin{equation}\label{eqn:leftInvProjConn}
\Omega = g_0^{-1} \, dg_0 + \Ad_{g_0}^{-1} C\left( h^{-1} \, dh \right)
\end{equation}
where $C : H \times G_0 \to \mathfrak{h}^* \otimes \mathfrak{g}$ is a function.
Check that $C$ is constant under the left $H$ action and under the right
$G_0$ action, so that $C$ is constant, an element of $\mathfrak{h}^* \otimes \mathfrak{g}$.
Moreover, this element, thought of as a linear map $\mathfrak{h} \to \mathfrak{g}$,
has image transverse to $\mathfrak{g}_0$.
Conversely, suppose that we pick any element $C \in \mathfrak{h}^* \otimes \mathfrak{g}$,
whose image is transverse to $\mathfrak{g}_0$.
We can construct a left invariant projective connection by equation~\ref{eqn:leftInvProjConn}.
It is easy to check that it is a projective connection.
\end{proof}
\begin{corollary}
For any $A \in \mathfrak{g}$, write $\bar{A}$ for the corresponding element of $\mathfrak{g}/\mathfrak{g}_0$.
A left invariant projective connection built from a linear map $C \in \mathfrak{h}^* \otimes \mathfrak{g}$
(with image of $C$ transverse to $\mathfrak{g}_0$) is normal just when
\[
\overline{C([A,B])} = \overline{\left[C(A),C(B)\right]},
\]
for any $A, B \in \mathfrak{g}$. It is flat just when $C$ is a Lie algebra homomorphism.
\end{corollary}
\begin{proof}
This is an easy calculation, given $\Omega$ in equation~\ref{eqn:leftInvProjConn}:
the curvature is 
\begin{align*}
d \Omega + \Omega \wedge \Omega &= \kappa \bar{\Omega} \wedge \bar{\Omega}
\\
&= 
\frac{1}{2}
\Ad_{g_0}^{-1} \left(
\left[
C\left(
h^{-1}  \, dh
\right),
C\left( 
h^{-1} \, dh
\right) 
\right]
-
C\left(
\left[
h^{-1} \, dh , 
h^{-1} \, dh
\right]
\right)
\right).
\end{align*}
so that 
\[
\kappa\left(\bar{A},\bar{B}\right)
=
\left[C(a),C(b)\right]-C\left(\left[a,b\right]\right)
\]
whenever $\bar{C}(a)=\bar{A}$ and $\bar{C}(b)=\bar{B}$,
for any $a,b \in \mathfrak{h}$.
\end{proof}
\begin{remark} The same approach will determine the isomorphism classes of left
invariant Cartan geometries of any type, and their curvature.
\end{remark}
We can always identify $\mathfrak{h}=\R{n}$, and then
we will have
\[
C(A) = 
\begin{pmatrix}
- \tr E(A) & F(A) \\
A & E(A)
\end{pmatrix},
\]
for some unique $E \in \mathfrak{h}^* \otimes \gl{\mathfrak{h}}$ and $F \in \mathfrak{h}^* \otimes \mathfrak{h}^*$.
This is normal just when
\[
\left(E(A) + \tr E(A)\right)B - \left(E(B) + \tr E(B)\right)A = \left[A,B\right].
\]
Its curvature is given by
\[
\kappa(A,B)=
\begin{pmatrix}
\kappa^0_0 & \kappa^0_{\bullet} \\
\kappa^{\bullet}_0 & \kappa^{\bullet}_{\bullet}
\end{pmatrix}
\]
where
\begin{align*}
\kappa^0_0 &= 
F(A)B - F(B)A + \tr E(\left[A,B\right]) \\
\kappa^0_{\bullet} &=
F(A) \left(E(B)+\tr E(B)\right)
-F(b) \left(E(A)+\tr E(A)\right)
\\
\kappa^{\bullet}_0 &=
\left(
E(A) + \tr E(A)
\right) B
-
\left(
E(B) + \tr E(B)
\right) A
-\left[A,B\right]
\\
\kappa^{\bullet}_{\bullet} &=
A \, F(B) - B \, F(A)
+ 
\left[E(A),E(B)\right]
-E\left(\left[A,B\right]\right).
\end{align*}
\subsection{Left invariant affine connections}
\begin{proposition}
The set of left invariant affine connections on a Lie group $H$
of dimension $n$ with Lie algebra $\mathfrak{h}$ is invariantly identified with $\mathfrak{h}^* \otimes \gl{\mathfrak{h}}$.
\end{proposition}
\begin{proof}
Take any left invariant connection on the tangent bundle of $H$, and identify
it as usual with a connection 1-form on $FH$. Let $\pi : FH \to H$ be the obvious
projection map. Define the soldering 1-form
$\omega \in \nForms{1}{FH} \otimes \R{n}$ by $v \hook \omega_u = u \left(\pi'(u) v\right)$.
It is easy to check that $r_g^* \omega = g^{-1} \omega$, for $g \in \GL{n,\R{}}$.
Then the connection 1-form $\gamma \in \nForms{1}{FH} \otimes \gl{n,\R{}}$
transforms in the adjoint representation under right $\GL{n,\R{}}$ action:
$r_g ^* \gamma = \Ad_g^{-1} \gamma$, and $d \omega + \gamma \wedge \omega = \frac{1}{2} T \omega \wedge \omega$
is the torsion of the connection.

Pick a point of $FH$, say $u_0$,
above the point $1 \in H$; so $u_0 : \mathfrak{h} \to \R{n}$. Fixing this identification,
we can say that $\omega \in \nForms{1}{FH} \otimes \mathfrak{h}$,
and that $\gamma \in \nForms{1}{FH} \otimes \gl{\mathfrak{h}}$.
The map $(h,g) in H \times \GL{\mathfrak{h}} \mapsto h_* g^{-1}u_0 \in FH$
is a diffeomorphism. So identify $FH$ with $H \times \GL{\mathfrak{h}}$ by
this diffeomorphism. This identifies $\omega_{(h,g)} = g^{-1} h^{-1} \, dh$,
clearly. Moreover, $\gamma=g^{-1} \, dg + \Ad_g^{-1} \Gamma\left(h^{-1} \, dh\right)$
for a unique constant choice of $\Gamma \in \mathfrak{h}^* \otimes \gl{\mathfrak{h}}$.
Moreover, reversing our steps ensures that all left invariant connections
on the tangent bundle occur in this manner, for a unique choice of $\Gamma$,
which can be selected arbitrarily.
\end{proof}
\begin{example}
For example, for any constant, we can take the choice 
$\Gamma = a \ad$, i.e. $\Gamma(A)=a \ad_A \in \gl{\mathfrak{h}}$,
giving a canonical choice of connection to the tangent bundle of any Lie
algebra. In particular, we can take $\Gamma=0$. Taking $\Gamma=0$ will give geodesic
flow $g(t)=g(0), h(t)=h(0) \, e^{t g(0) A}$, for any constant $A \in \mathfrak{h}$,
so a complete connection.
\end{example}
\begin{example}
If we pick any nondegenerate quadratic form on the Lie algebra $\mathfrak{h}$,
then we can define $\ad^t$ by 
\[
\left<\ad^t_A B, C\right> = \left<B, \ad_A C\right>
\] 
for $A, B, C \in \mathfrak{h}$, and define $\ad'$ by
\[
\ad'_A B= \ad^t_{B} A.
\]
Then the Levi--Civita connection of the induced left invariant Riemannian
metric on $H$ is given by $\Gamma=\frac{1}{2}\left(\ad-\ad^t-\ad'\right)$.

If the metric is bi-invariant, then $\ad^t=-\ad$ and $\ad'=\ad$, so one
has $\Gamma=\frac{1}{2} \ad$.
\end{example}
\begin{corollary} A left invariant connection on the tangent bundle of a Lie group
is torsion-free just when $\Gamma=\frac{1}{2}\ad + S$, where
$S \in \Sym{2}{\mathfrak{h}}^* \otimes \mathfrak{h}$.
\end{corollary}
\begin{proof}
One easily calculates $d \omega + \gamma \wedge \omega$,
and finds that $\Gamma$ gives a torsion-free connection
just when $\Gamma(A)B-\Gamma(B)A-[A,B]$ is symmetric
in $A$ and $B$. Indeed the torsion of the connection is
\[
\frac{1}{2} T \omega \wedge \omega = 
g^{-1} \left( 
\Gamma\left(h^{-1} \, dh\right) \wedge h^{-1} \, dh  - \frac{1}{2} \left[h^{-1} \,  dh, h^{-1} \, dh \right] 
\right)
\]
so that
\[
T\left(A,B\right) = \Gamma(A)B-\Gamma(B)A - \left[A,B\right].
\]
\end{proof}
\begin{corollary} A torsion-free left invariant connection corresponds
to a choice of $S \in \Sym{2}{\mathfrak{h}}^* \otimes \mathfrak{h}$
invariant under the adjoint action.
\end{corollary}
\begin{corollary} A left invariant connection on the tangent bundle of a Lie group
is flat just when $\Gamma$ is a Lie algebra morphism.
\end{corollary}
\begin{proof}
It is easy to compute that the curvature is
\begin{align*}
d \gamma + \frac{1}{2} \left[\gamma,\gamma\right] &=
\frac{1}{2} \kappa \omega \wedge \omega \\
&=
\frac{1}{2} 
\Ad_{g}^{-1} 
\left(
\left[
\Gamma\left(h^{-1} \, dh \right),
\Gamma\left(h^{-1} \, dh \right)
\right]
-
\Gamma
\left(
\left[h^{-1} \, dh, h^{-1} \, dh \right]
\right)
\right),
\end{align*}
so that 
\[
\kappa(A,B) = 
\left[
\Gamma\left(A \right),
\Gamma\left(B\right)
\right]
-
\Gamma
\left(
\left[A,B \right]
\right).
\]
\end{proof}
\begin{proposition}
If $\Gamma \in \mathfrak{h}^* \otimes \gl{\mathfrak{h}}$ determines 
a left invariant affine connection on a Lie group $H$ of dimension $n$ with Lie 
algebra $\mathfrak{h}$, then 
\[
C(A) = 
\begin{pmatrix}
-\frac{1}{n+1} \tr \Gamma(A) & \frac{1}{2(n+1)} \tr \Gamma \circ \ad_A \\
A & \Gamma(A)-\frac{1}{n+1} \tr \Gamma(A)
\end{pmatrix}
\]
determines the corresponding projective connection, which is normal
just when $\Gamma$ is torsion-free.
\end{proposition}
\begin{proof}
One easily calculates out the geodesic equation, to see that it agrees,
and checks that the curvature is suitably trace-free.
\end{proof}
\begin{lemma}  Let $\Killing$ be the Killing form of a Lie algebra $\mathfrak{h}$:
\[
\Killing(A,B)=\tr \ad_A \ad_B.
\]
The symmetrized Ricci curvature of the natural torsion-free connection given by $\Gamma = \frac{1}{2} \ad$
(which is the Levi--Civita connection of the Killing form metric on any semisimple
Lie group) is $r=-\frac{1}{4} \Killing$.
\end{lemma}
\begin{proof}
Pick a basis of $\mathfrak{h}$, say $e_1, e_2, \dots, e_n$.
Suppose that the structure constants are $c^i_{jk}$, so that $\left[e_i, e_j\right] = c^k_{ij} e_k$.
Then $\Killing_{ij} = c^k{i \ell} c^{\ell}_{jk}$. The curvature is
\begin{align*}
K^i_{jk \ell} &= 
\left[
\frac{1}{2} \ad_{e_k} , \frac{1}{2} \ad_{e_{\ell}}
\right]
-
\frac{1}{2}
\ad_{
\left[
e_k, e_{\ell}
\right]
}
\\
&=
\frac{1}{4} c^i_{k m} c^m_{\ell j}
-
\frac{1}{4} c^i_{\ell m} c^m_{k j}
-
\frac{1}{2} c^m_{k \ell} c^i_{m j}.
\end{align*}
Therefore the symmetrized Ricci tensor is
\begin{align*}
K_{j \ell} &=
\frac{1}{2}\left(K^i_{j i  \ell} + K^i_{\ell i j}\right)
\\
&=
\frac{1}{2}\left(
\frac{1}{4} c^i_{im} c^m_{\ell j}
-
\frac{1}{4} c^m_{\ell i} c^i_{jm}
+
\frac{1}{4} c^i_{im} c^m_{j \ell}
-
\frac{1}{4} c^m_{j i} c^i_{\ell m}
\right) \\
&=
-\frac{1}{4} c^m_{j i} c^i_{\ell m}
\\
&= - \frac{1}{4} \Killing_{j \ell}.
\end{align*}
\end{proof}
\begin{example} Consider this same natural torsion-free connection given by $\Gamma = \frac{1}{2} \ad$.
Its geodesics are precisely the one-parameter subgroups (an easy calculation). Its symmetrized 
Ricci curvature is the Killing form, which is bi-invariant. Therefore the Ricci
curvature in the direction of a given geodesic is constant. This ensures
that this torsion-free connection is projectively complete just precisely
when the Killing form is positive definite, i.e. precisely on the compact
semisimple Lie groups. 
\end{example}
\begin{example}
For example, $\SL{2,\R{}}$ has projectively complete geodesics
given precisely by the subgroups conjugate to $\SO{2}$, and has hyperbolic
geodesics in the tangent directions of the 2-dimensional subgroups,
and parabolic geodesics precisely in the directions of nilpotent elements, 
i.e. in directions conjugate to the subgroup generated by
\[
\begin{pmatrix}
0 & 1 \\
0 & 0 
\end{pmatrix}.
\]
There is no other bi-invariant torsion-free connection on $\SL{2,\R{}}$,
because (by Clebsch--Gordan) there are no nonzero
elements of $\Sym{2}{\slLie{2,\R{}}}^* \otimes \slLie{2,\R{}}$
fixed under $\SL{2,\R{}}$.
\end{example}
\begin{example}
Looking at the representations of $\slLie{2,\R{}}$, we see that the biinvariant
projective connections on $\SL{2,\R{}}$ are precisely those of the form
\[
C(A)
=
\begin{pmatrix}
0 & q \, A^* \\
A & p \, \ad_A
\end{pmatrix}
\]
where $A^*(B)=\Killing(A,B)$ is the dual covector in the Killing form,
and $p$ and $q$ are arbitrary constants. The curvature is
\[
\kappa(A,B)
=
\begin{pmatrix}
0 & 2pq \, \left[A,B\right]^* \\
(2p-1) \, \left[A,B\right] & q \left(A \otimes B^*  - B \otimes A^* \right)
+ p(p-1) \ad_{[A,B]}
\end{pmatrix}.
\]
This projective connection is torsion-free just when $p=\frac{1}{2}$. Its
symmetrized Ricci curvature is
\[
r(A,B)=\left(
(n-1)q + p(p-1)  
\right) \Killing(A,B).
\]
However, our theorem on Ricci curvature only applies when the
projective connection is normal, and this happens only
for the example we have already calculated. To see what is
at issue more clearly, lets write the equations of the geodesic flow.
First, for $g_0 \in G_0$, write
\[
g_0 = 
\begin{pmatrix}
a & b^* \\
0 & c
\end{pmatrix}.
\]
The equations of geodesic flow, expressed 
in $\Omega = g_0^{-1} \, dg_0 + \Ad_{g_0}^{-1} C\left( h^{-1} \, dh \right)$,
become
\begin{align*}
da &= \left<b,A\right>  \\
db &= - \frac{q}{\det c} c^t c A \\
dc &= - \frac{1}{\det c} \left(
\left(cA\right) \otimes b + p \, \ad_{cA} c 
\right) \\ 
h^{-1} \, dh &= \frac{1}{\det c} cA,
\end{align*}
where $A$ is any fixed element of $\slLie{2,\R{}}$.
We can see that solving for 
$b$ and $c$ first should enable us to solve the
whole system. Nonetheless, the author cannot see
how to solve these ordinary differential equations,
or how to estimate the time during which the solutions
remain defined.
\end{example}
\begin{remark} Blumenthal \cite{Blumenthal:1987} proves that
totally geodesic fiber bundles with projectively complete
total space have projectively complete base space. For example,
the Hopf fibration ensures projective completeness of complex
projective space. But complex projective space has positive 
Ricci curvature, so we don't actually need Blumenthal's result
to see projective completeness.
\end{remark}
\section{Jacobi vector fields}
Consider a family of geodesics in a manifold $M$ with projective
connection. Let $S$ be a surface, with a submersion
$t : S \to \R{}$ whose fibers are curves $S_t \subset S$.
Map $S \to M$ so that each fiber $S_t$ maps into a geodesic. 
Let $E_S$ be the union of the bundles $E_{S_t}$.
On $E_S$,
\begin{align*}
\omega^I &= a^I \, dt \\
\gamma^I_1 &= a^I_1 \, dt.
\end{align*}
Take the exterior derivative to find
\begin{align*}
\bar{\nabla}
\begin{pmatrix}
a^I \\
a^I_1
\end{pmatrix}
&=
d
\begin{pmatrix}
a^I \\
a^I_1
\end{pmatrix}
+
\begin{pmatrix}
\gamma^I_J & 0 \\
- \delta^I_J \omega_1 & \gamma^I_J - \delta^I_J \gamma^1_1
\end{pmatrix}
\\
&=
\begin{pmatrix}
K^I_{1J} & \delta^I_J \\
K^I_{11J} & 0
\end{pmatrix}
\begin{pmatrix}
a^J \\ 
a^J_1
\end{pmatrix}
\omega^1
+
\begin{pmatrix}
\bar{\nabla}_t a^I \\
\bar{\nabla}_t a^I_1
\end{pmatrix}
\, dt.
\end{align*}
where $\bar{\nabla}$ represents the covariant derivative
in $E_S$ calculated in the coframing $dt,\omega^1,\omega_1,\dots$.
Therefore each motion $S$ through geodesics determines,
above each geodesic $S_0$ some functions 
$\left(a^I,a^I_1\right) : E_{S_0} \to \R{2(n-1)}$, 
satisfying the equations
\begin{equation}\label{eqn:JVFone}
d
\begin{pmatrix}
a^I \\
a^I_1
\end{pmatrix}
+
\begin{pmatrix}
\gamma^I_J & 0 \\
-\delta^I_J \omega_1 & \gamma^I_J - \delta^I_J \gamma^1_1
\end{pmatrix}
\begin{pmatrix}
a^J \\
a^J_1
\end{pmatrix}
=
\begin{pmatrix}
K^I_{1J} & \delta^I_J \\
K^I_{11J} & 0
\end{pmatrix}
\begin{pmatrix}
a^J \\ 
a^J_1
\end{pmatrix}
\omega^1.
\end{equation}
Moreover, under the action of $G_+$,
\begin{equation}\label{eqn:JVFtwo}
r_g^* 
\begin{pmatrix}
a^I \\
a^I_1
\end{pmatrix}
=
\begin{pmatrix}
\left(g^{-1}\right)^I_J & 0 \\
\left(g^{-1}\right)^I_J g^0_1 & \left(g^{-1}\right)^I_J g^1_1 \\
\end{pmatrix}
\begin{pmatrix}
a^J \\
a^J_1
\end{pmatrix}.
\end{equation}
\begin{definition}
If $C \subset M$ is a geodesic,  
and $\left(a^I, a^I_1\right) : E_C \to \R{2(n-1)}$
is a solution to equations~(\ref{eqn:JVFone}),~(\ref{eqn:JVFtwo}),
call $\left(a^I, a^I_1\right)$ a \emph{Jacobi vector field} (even
if it doesn't arise from any actual
variation through geodesics).
\end{definition}
\begin{lemma}
Suppose that a connected manifold $M$ has
a projective connection $E \to M$.
Any infinitesimal symmetry of a projective 
connection $E \to M$ determines
and is determined by a unique vector
field on $M$.
\end{lemma}
\begin{proof}
The equations of an infinitesimal symmetry
are $\LieDer_X \Omega=0$ and $r_{g*} X = X$,
$X$ a vector field on $E$. But then
if we let $X^{\bullet} = X \hook \omega^{\bullet}$,
we find immediately that $X^{\bullet}$ represents
a vector field on $M$. Moreover, if $X^{\bullet}=0$,
then $X \hook \Omega : E \to \mathfrak{g}$
is a $G_0$-equivariant function:
\[
r_g^* \left( X \hook \Omega \right ) = \Ad_g^{-1} \left( X \hook \Omega \right),
\]
but is also invariant under the flow of $X$.
Moreover, for any vector $\vec{A} \in \mathfrak{g}$,
\[
\LieDer_{\vec{A}} X \hook \Omega = - \Ad_A X \hook \Omega.
\]
Pick any point $e \in E$,
and let $A=X(e) \hook \Omega(e)$.
Consider the vector field $\vec{A}$.
It agrees with $X$ at $e$, and
also satisfies 
$r_g^* \left ( \vec{A} \hook \Omega \right ) = \Ad_g^{-1} \left ( \vec{A} \hook \Omega \right)$,
so $\vec{A}$ agrees with $X$ up the fiber of $e$.
Moreover, $\vec{A}$ and $X$ have the same
brackets with $\vec{B}$ for any $B \in \slLie{n+1}$.
Therefore they agree above the path component of 
$M$ containing $\pi(e)$, and since $M$ is connected
we find $X=\vec{A}$. 
But $X$ is $G_0$-invariant, so $A$ must be
$G_0$-invariant, i.e $A \in \mathfrak{g}_0$
belongs to the center of $\mathfrak{g}_0$.
Check that the center is $0$.
\end{proof}
Along a curve $S_0$, we can construct the normal 
bundle $\nu S_0 = \left.TM\right|_{S_0}/TS_0$.
Given any section $a$ of the normal bundle,
define functions $a^I$ by 
\[
a^I = a \hook \underline{\omega}^I,
\]
which is well defined because $\underline{\omega}^I$
vanishes on $TS_0$.
\begin{lemma}
Given an immersed curve $S_0 \subset M$, 
and a section $a$ of its normal bundle,
the functions $a^I$ satisfy
\[
r_g^* a^I = \left(g^{-1}\right)^I_J a^J
\]
for $g \in G_+$.
Conversely, given functions satisfying
these equations, they determine
a section of the normal bundle.
\end{lemma}
\begin{corollary}
A Jacobi vector field determines and is determined by a section of
the normal bundle. For a family of geodesics $x : S^1 \times \R{} \to M$,
this section is
\[
a(x)=
\left.\pd{x}{t}\right|_{t=0} \mod x_* TS^1,
\]
the normal component of the velocity.
\end{corollary}
\begin{proof}
At each point $e \in E_S$, we see that 
\begin{align*}
\pd{}{t} \hook \underline{\omega}^I
&=
\pd{}{t} \hook a^I \, dt \\
&= a^I.
\end{align*}
\end{proof}
Let $N \subset G$ be the subgroup
acting trivially on the projective line
$\Proj{1}$ containing $\left[e_0\right]$
and $\left[e_1\right]$,
and $N_+ = N \cap G_+$. 
Let $\bar{G}=G/N$ and $\bar{G}_+=G_+/N_+$.

Given an immersed curve $C$, the quotient $\bar{E}_C = E_C/N_+ \to C$
is a principal right $\bar{G}_+$-bundle, and when
equipped with the 1-form 
\[
\bar{\Omega} = \Omega \mod \mathfrak{n}_+ =
\begin{pmatrix}
-\frac{1}{2} \gamma^1_1 & \omega_1 \\
\omega^1 & \frac{1}{2} \gamma^1_1
\end{pmatrix}
\]
is a flat projective connection on $C$.

The sheaf of infinitesimal symmetries 
on $\tilde{C}$ is the sheaf of solutions of 
a system of linear ordinary differential equations,
so by Picard's theorem the local solutions
extend globally, never becoming multivalued
because $\tilde{C}$ is simply connected.
\begin{lemma}
$C$ is complete just when every infinitesimal
symmetry of the projective connection
on $\tilde{C}$ is complete.
\end{lemma}
\begin{proof}
Follows from the classification of projective
connections on curves.
\end{proof}
\begin{lemma}
If $C$ is complete then the map 
$\tilde{C} \to \Proj{1}$ preserves
and reflects infinitesimal symmetries.
\end{lemma}
\begin{proof}
This is clear from the classification.
Alternatively: the infinitesimal symmetries pullback to 
infinitesimal symmetries, because 
$\tilde{C} \to \Proj{1}$ is a covering map.
Infinitesimal symmetries comprise a 3 dimensional
vector space, so they must match precisely. 
\end{proof}
The equations of an infinitesimal symmetry
\[
\LieDer_X 
\begin{pmatrix}
\omega^1 \\
\gamma^1_1 \\
\omega_1
\end{pmatrix},
\]
if we set 
\[
\begin{pmatrix}
X^1 \\
X^1_1 \\
X_1
\end{pmatrix}
=
X \hook
\begin{pmatrix}
\omega^1 \\
\gamma^1_1 \\
\omega_1
\end{pmatrix}
\]
give us
\begin{equation}\label{eqns:InfSymm}
d
\begin{pmatrix}
X^1 \\
X^1_1 \\
X_1
\end{pmatrix}
=
\begin{pmatrix}
X^1_1 \omega^1 - X^1 \gamma^1_1 \\
2 X^1 \omega_1 - 2 X_1 \omega^1 \\
X_1 \gamma^1_1 - X^1_1 \omega_1
\end{pmatrix}.
\end{equation}
\begin{lemma}
If $X$ is an infinitesimal symmetry,
let 
\[
\mathbb{B}(X)=2X_1 X^1 + \frac{1}{2}\left(X^1_1\right)^2.
\]
Then $\mathbb{B}(X)$ is a constant.
\end{lemma}
\begin{proof}
This is just the Killing form applied to
$X \hook \Omega'$. It is also easy to
check by calculating the exterior derivative of $\mathbb{B}(X)$.
\end{proof}
\begin{definition}
An infinitesimal symmetry $X$ is called
\emph{elliptic} if $\mathbb{B}(X)<0$, \emph{parabolic}
if $\mathbb{B}(X)=0$ and $X$ is not everywhere $0$, 
and \emph{hyperbolic} if $\mathbb{B}(X)>0$.
\end{definition}
\begin{lemma}
Every elliptic infinitesimal symmetry has
no zeros. Zeros of a hyperbolic infinitesimal
symmetry (if there are any zeros) are of
order 1. Zeros of a parabolic infinitesimal
symmetry (if there are any zeros) are of
order precisely 2.
\end{lemma}
\begin{proof}
Zeros here mean on $\tilde{C}$, so
equate upstairs on $\bar{E}_{\tilde{C}}$ to
zeros of $X^1$.

Either use the fact that this lemma holds true
on $\Proj{1}$, and local isomorphism
of sheaves of infinitesimal symmetries,
or more simply note that at a zero:
\[
\mathbb{B}(X)=\frac{1}{2}\left(X^1_1\right)^2,
\]
which can't be negative, ruling out ellipticity.
Moreover, it is positive (hyperbolicity) just
when $X^1_1 \ne 0$, i.e. just when $dX^1 \ne 0$,
a zero of order 1. Any zero of order higher than
2 would ensure that the differential equations
(\ref{eqns:InfSymm}) for infinitesimal symmetries 
have the same initial conditions as the 0 infinitesimal
symmetry.
\end{proof}
\begin{lemma}
$C$ is complete just when
all parabolic infinitesimal
symmetries of $\tilde{C}$ are complete.
\end{lemma}
\begin{proof}
The infinitesimal symmetries form a Lie
algebra spanned by the parabolic
ones. By Palais' theorem \cite{Palais:1957},
a finite dimensional Lie algebra generated
by complete vector fields consists entirely
of complete vector fields.

Alternatively, just look again at the
classification of projective connections
on curves.
\end{proof}
\section{Normal projective connections on surfaces}
\begin{definition}
A projective connection is called \emph{normal}
if
\begin{align*}
0 &= K^i_{jk} \\
0 &= K^i_{jil} \\
\intertext{and}
0 &= K^i_{jkl} + K^i_{ljk} + K^i_{klj}.
\end{align*}
\end{definition}
\begin{lemma}[Cartan \cite{Cartan:1992}]
Given a projective connection, there is
a unique normal projective connection
with the same unparameterized geodesics.
\end{lemma}
\begin{lemma}
Let $M$ be a surface with normal projective connection
$E \to M$, and  $C$ a periodic geodesic.
If $C$ has a nonzero Jacobi
vector field, which vanishes at some point of $C$,
then $C$ is a complete geodesic.
\end{lemma}
\begin{proof}
The equations of a Jacobi vector field
are identical to the flat case, since
the relevant curvature vanishes. Therefore
under development the Jacobi vector fields
are identified locally with Jacobi vector fields
on the model. 
Any parabolic infinitesimal symmetry on $\ucProj{1}$ has 
a zero between any two zeros of a
Jacobi vector field. But
$\tilde{C} \subset \ucProj{1}$
is just an open interval, with the same differential
system for Jacobi vector fields, so the same is
true on $\tilde{C}$. But this forces every parabolic
infinitesimal symmetry to have zeros arbitrarily
far along $\tilde{C}$ in both directions, since
the Jacobi vector fields have zeros in $C$,
so periodically placed zeros in $\tilde{C}$
arbitrarily far along.
This forces the parabolic infinitesimal
symmetries to be complete, since the
flow of a parabolic infinitesimal symmetry
will drive us toward its next zero.
\end{proof}
\section{Tameness}
We follow LeBrun \& Mason \cite{LeBrunMason:2002} closely here;
keep in mind that their paper treats projective
structures, rather than the more general concept
of projective connection, so one has to check
that the results quoted below hold, with the same
proofs, for arbitrary projective connections. Our aim is to show that 
every normal projective connection on a surface
whose geodesics are all closed has a nonzero
Jacobi vector field on each geodesic, i.e. lots of
motions through closed geodesics.

Given a projective connection $\pi : E \to M$,
define a map $\Pi : E \to \Proj{}TM$ by
requiring that for any geodesic
$C$, $e \in E_C \mapsto \Pi(e)= T_{\pi(e)} C \in \Proj{}TM$.
\begin{lemma}
This map $\Pi$ is well-defined and smooth and 
identifies $\Proj{}TM=E/G_+$.
The foliation of $E$ by the flow lines
of the geodesic flow descends to a 
foliation of $\Proj{}TM$, whose
leaves project to the geodesics in 
$M$.
\end{lemma}
\begin{proof}
To show that $\Pi$ is well-defined,
we have only to show that for 
each point $e \in E$, there is 
a geodesic $C \subset M$ with
$e \in E_C$. But we can just take
the integral manifold of the geodesic
exterior differential system in $E$
passing through $e$, and it will
project an appropriate geodesic $C$.
Indeed the
geodesic flow line through $e$
projects to an appropriate geodesic. This
makes clear the smoothness of $\Pi$.

Suppose that we have two points $e_0,e_1 \in E$
with $\Pi\left(e_0\right)=\Pi\left(e_1\right)=\ell \in \Proj{}TM$.
Then the integral manifolds $E_{C_0}$ and $E_{C_1}$ 
of the exterior differential system for geodesics
with $e_j \in E_{C_j}$ must project to tangent
geodesics: $T_{c_0} C_0=T_{c_1} C_1$ for some
$c_0$ and $c_1$. The 1-forms $\underline{\omega}^I\left(e_j\right)$
must therefore be linear multiples of one another,
which forces $e_0$ and $e_1$ to be in the same
$G_+$-orbits, by looking at how the $\omega^i$
transform under right $G_+$-action. Therefore
$E_{C_0}=E_{C_1}$. So $\Pi : E/G_+ \to \Proj{}TM$
is 1-1. Under right $G_0$-action,
$r_{g_0} \vec{A} = \overrightarrow{\Ad_{g_0}^{-1} A}$
ensuring that $\Pi$ is onto, since this
action acts transitively on semibasic directions.
To ensure that the inverse map $\Proj{}TM \to E/G_+$
is smooth, take any local section of 
$E \to E/G_+$, and attach to each point
$\ell \in \Proj{}TM$ the associated point of $E$;
this is the point satisfying the equations
$\bar{\omega}^I\left(e\right)=0$ on $\ell$.
Once again, examining the right action of $G_0$,
it is easy to check that this point $e$ is uniquely
determined and smoothly so.
\end{proof}
\begin{definition}
A projective connection $E \to M$ 
is called \emph{tame} if the foliation
of $\Proj{}TM$ is locally trivial.  
It is called \emph{Zoll} if all geodesics
are embedded closed curves.
\end{definition}
\begin{lemma}[LeBrun \& Mason \cite{LeBrunMason:2002}]
Any Zoll projective connection on a compact
surface is tame. The only surfaces
which admit Zoll projective connections
are the sphere and the real projective plane.
\end{lemma}
\begin{lemma}[LeBrun \& Mason \cite{LeBrunMason:2002}]
Let $M$ be a surface bearing a Zoll
projective connection.
The space of unoriented connected geodesics $\Lambda$
(i.e. the space of leaves of the foliation
of $\Proj{}TM$) is diffeomorphic to $\Proj{2}$.
The map $\Proj{}TM \to \Lambda$ 
is a smooth fiber bundle, taking
$T_c C \mapsto C$ (for $C$ any geodesic and $c \in C$).
\end{lemma}
\begin{lemma}
Let $\Lambda$ be the space of geodesics of 
a tame Zoll projective connection $M$.
The map $\Proj{}TM \to \Lambda$ linearly identifies
Jacobi vector fields with tangent vectors to $\Lambda$.
\end{lemma}
\begin{proof}
Every vector in $T \Lambda$ gives rise 
to an infinitesimal motion through geodesics,
with nonvanishing normal velocity, and
therefore by dimension count must account
for all of the Jacobi vector fields,
since the ordinary differential equation
for Jacobi vector fields has well defined
initial value problem.
\end{proof}

\begin{corollary}
Zoll normal projective connections on surfaces are complete.
\end{corollary}

\begin{example}
Zoll \cite{Zoll:1903} provides the following examples
of Zoll metrics: for any odd function $f : [-1,1] \to (-1,1)$,
with $f(z)=-f(-z)$ and $f(-1)=f(1)=0$,
\[
ds^2 = \frac{\left(1+f(z)\right)^2}{1-z^2} \, dz^2 + 
\left(1-z^2\right) \, d\theta^2
\]
is a Zoll metric on the 2-sphere (i.e. all
geodesics are embedded and periodic), 
for $(z,\theta) \in [-1,1] \times [0,2 \pi]$
longitude and latitude coordinates.
The curvature is
\[
\kappa =
\frac{f(z)+1-zf'(z)}{\left(f(z)+1\right)^3}.
\]
For instance, if $f(z) \equiv 0$, then $\kappa \equiv 1$
giving the standard metric on the sphere.
On the other hand, if $f(z)$ has 
large first derivative and small value 
at some point $z_0$ away from $0$ (and, being
odd, has the same behaviour at $-z_0$), 
then we will find negative curvature near $z=z_0$.
For example
\[
f(z)=\cos\left(\frac{\pi}{2} z\right) e^{- \alpha z^2} 
\sin\left(2 \arctan\left(\beta \left(z-z_0\right)\right)
+2 \arctan\left(\beta \left(z+z_0\right)\right)
\right),
\]
with large constants $\alpha$ and $\beta$,
so that $\beta e^{-\alpha z_0^2}$ is large,
and any $z_0 \ne 0$ with $-1 < z_0 < 1$, 
gives curvature close to 1 near $z= \pm 1$,
and negative near $z=\pm z_0$, on a rotationally
symmetric surface which is as close as we like
to the usual metric on the sphere near the poles.
Indeed $\alpha=1, \beta=1/4, z_0=1/2$ already provides
two bands of negative curvature. Similarly, 
we can construct any even number of negative
curvature bands, at prescribed locations.
\end{example}
\begin{remark}
LeBrun \& Mason \cite{LeBrunMason:2002}
provide explicit examples of Zoll affine connections
on the sphere, which must also 
determine complete projective connections by our results.
\end{remark}
\begin{remark}
It is not clear if there is any simple relationship between
projective and affine completeness. It is not known whether
torsion-free Zoll affine connections are complete.
\end{remark}
\section{Conclusions}
There is no clear common thread between
our examples. Perhaps there is always
a positive Ricci curvature affine connection
for any Zoll normal projective connection.
Perhaps Zoll normal projective connections
are always complete.
It seems a reasonable conjecture that any projective connection can be
perturbed to a projective
connection all of whose geodesics are hyperbolic,
without altering its unparameterized geodesics.
It appears unlikely that a projective connection can be 
made complete by such a perturbation.
\bibliographystyle{amsplain} 
\bibliography{ZollComplete}
\end{document}